\documentclass[reqno,10.5pt]{amsart}

\usepackage{amsthm,amsmath,amssymb}
\usepackage{mathrsfs,amsfonts,dsfont,functan,extarrows,mathtools}
\usepackage[colorlinks]{hyperref}

\usepackage{marginnote}
\usepackage{xcolor}

\makeatletter

\newcommand{\Rmnum}[1]{\expandafter\@slowromancap\romannumeral #1@}
\makeatother

\def\Xint#1{\mathchoice	{\XXint\displaystyle\textstyle{#1}}%	
{\XXint\textstyle\scriptstyle{#1}}%	
{\XXint\scriptstyle\scriptscriptstyle{#1}}%	
{\XXint\scriptscriptstyle\scriptscriptstyle{#1}}%	
\!\int}\def\XXint#1#2#3{{\setbox0=\hbox{$#1{#2#3}{\int}$}	
\vcenter{\hbox{$#2#3$}}\kern-.5\wd0}}\def\dashint{\Xint-}

\newtheorem{theorem}{Theorem}[section]
\newtheorem{proposition}[theorem]{Proposition}

\newtheorem{lemma}[theorem]{Lemma}

\numberwithin{equation}{section}
\allowdisplaybreaks

\newcounter{wronumber}

\begin{document}

\title[]{Existence of the $(\alpha,\beta)$-Ricci-Yamabe flow on closed manifolds}

\author[L. Zhang]{Liangdi Zhang}
\address[Liangdi Zhang]{\newline 1. Yanqi Lake Beijing Institute of Mathematical Sciences and Applications, Beijing 101408, P. R. China
\newline 2. Yau Mathematical Sciences Center, Tsinghua University, Beijing 100084, P. R. China}
\email{ldzhang91@163.com}

%\thanks{First Draft}
%\date{December 2021}

%%%%%%%%%%%%%%%%%%%%%%%%%%%%%%%%%%%%%%%%%%%%%%%%%%%%%%%%%%%%%%%%%%%%%%%%%%%%%%
%%%%%%%%%%%%%%%%%%%%%%%%%%%%%%%%%%%%%%%%%%%%%%%%%%%%%%%%%%%%%%%%%%%%%%%%%%%%%%

\begin{abstract}
On a smooth closed Riemannian manifold, we show short time existence of smooth solutions to the $(\alpha,\beta)$-Ricci-Yamabe flow, which is a natural generalization of the Ricci flow and the Yamabe flow. We also establish some long time existence theorems for the closed $(\alpha,\beta)$-Ricci-Yamabe flow by estimating its curvatures.
\vspace*{5pt}


\end{abstract}

\maketitle

\tableofcontents

%%%%%%%%%%%%%%%%%%%%%%%%%%%%%%%%%%%%%%%%%%%%%%%%%%%%%%%%%%%%%%%%%%%%%%%%%%%%%%
%%%%%%%%%%%%%%%%%%%%%%%%%%%%%%%%%%%%%%%%%%%%%%%%%%%%%%%%%%%%%%%%%%%%%%%%%%%%%%
\section{Introduction} % (fold)
\label{sec:1}
Let $(M^n,g(x,t))$ be an $n$-dimensional smooth closed Riemannian manifold with the metric $g(x,t)$ evolves under the $(\alpha,\beta)$-Ricci-Yamabe flow
\begin{equation}\label{1.1}
%\begin{cases}
\frac{\partial}{\partial t}g=-2\alpha Ric-\beta Rg,
%g(0)=g_0,
%\end{cases}
\end{equation}
where $\alpha$ and $\beta$ are real constants.

The $(\alpha,\beta)$-Ricci-Yamabe flow was introduced by G\"uler and Grasmareanu \cite{RYmap} in 2019. This family of geometric flows includes or relates to two famous geometric flows: the Ricci flow ($\alpha=1$ and $\beta=0$) and the Yamabe flow ($\alpha=0$ and $\beta=1$).

Using the Nash-Moser inverse function theorem, Hamilton \cite{3D} proved short time existence for smooth solutions to the Ricci flow
\[\frac{\partial}{\partial t}g=-2 Ric\]
on a closed manifold for the first time. Shortly after that, DeTurck \cite{DT83} subsequently simplified the short time existence proof by modifying the flow and showing that the Ricci flow could be replaced by an equivalent PDE which is strictly parabolic. More generally, Shi \cite{Shi89} proved short time existence for the Ricci flow on a complete noncompact Riemannian manifold if the initial Riemannian curvature is bounded.

In 1995, Hamilton \cite{Ha95} provided a long time existence criterion that a closed Ricci flow defined on a maximal interval $[0,T)$ satisfies either $T=\infty$ or the maximum of the norm of the Riemannian curvature blows up at the finite time $T$. About a decade later, \v{S}e\v{s}um \cite{Se05} proved that if the Ricci curvature is uniformly bounded  along the closed Ricci flow defined on a finite time interval, then the Riemannian curvature stays uniformly bounded along the flow by using Perelman's noncollapsing theorem (see \cite{pere02}) and Hamilton's compactness theorem (see \cite{Hacpt}). \v{S}e\v{s}um's result still holds on the noncompact setting that the uniformly bounded Ricci curvature implies the uniformly bounded Riemannian curvature along the Ricci flow on a complete noncompact manifold with bounded Riemannian curvature at the initial time (see the work of Ma-Cheng \cite{MC11}, Kotschwar-Munteanu-Wang \cite{KMW16}, and Hsu \cite{Hsu19}). 

%In 2019,  provided another proof by using the De Giorgi iteration method.

%There are also long time existence theorems for the Ricci flow with integral curvature conditions, for example, \cite{wb08,MC11,dimatteo}.
There is a well-known conjecture on a closed Ricci flow $g(\cdot,t)$ $0\leq t<T<\infty$ that a uniform bound for the scalar curvature is enough to extend Ricci flow over time $T$. This conjecture is achieved in dimension 3 by Hamilton-Ivey's pinching estimate (see \cite{ha99,ivey93}). Zhang \cite{zz10} and Enders-M\"uller-Topping \cite{emt11} partially settled this conjecture for K\"ahler-Ricci flow and Type \Rmnum{1} Ricci flow, respectively. Wang \cite{wb08} proved that the closed Ricci flow $g(\cdot,t)$ $0\leq t<T<\infty$ can be smoothly extended past $T$ if the Ricci curvature tensor is uniformly lower bounded and the scalar curvature satisfies some space-time integral bounds. This result was generalized to the complete case by Di Matteo \cite{dimatteo} with additional conditions in 2021.

Ma and Cheng \cite{MC11} derived a smooth extension result for a closed Ricci flow with finite $L^\frac{n+2}{2}$ norms of the scalar curvature and Weyl tensor, and showed the Riemannian curvature of a complete Ricci flow stays uniformly bounded along the flow if the sectional curvature is bounded at the initial time and the scalar curvature and Weyl tensor are uniformly bounded. For a closed Ricci flow $g(\cdot,t)$ defined on a maximal interval $[0,T)$,  Cao \cite{cao11} proved either the scalar curvature blows up at the finite time $T$ or
\[\limsup_{t\rightarrow T}R<\infty\ \text{but}\ \limsup_{t\rightarrow T}\frac{|W|}{R}=\infty,\]
where $W=\{W_{ijkl}\}$ is the Weyl tensor of $g(\cdot,t)$ defined by
\begin{eqnarray*}
W_{ijkl}&=&R_{ijkl}-\frac{1}{n-2}(g_{ik}R_{jl}-g_{il}R_{jk}-g_{jk}R_{il}+g_{jl}R_{ik})\notag\\
&&+\frac{1}{(n-1)(n-2)}R(g_{ik}g_{jl}-g_{il}g_{jk}).
\end{eqnarray*}
In 2018, Li \cite{lirhf18} extended Cao's result \cite{cao11} to the Ricci-harmonic flow
\[
\begin{cases}
\frac{\partial}{\partial t}g=-2Ric+2a(t)\nabla\phi\otimes\nabla\phi,\\
\frac{\partial}{\partial t}\phi=\Delta\phi,
\end{cases}
\]
where $a(t)$ is a time-dependent positive constant and $\phi$ is a family of functions on $M\times[0,T)$ (with $0<T\leq\infty$). For more related works on the existence of Ricci-harmonic flow, please refer to \cite{list08,muller12,chengzhu13,LY18,LZ21ax,wuzheng20}.

The Yamabe flow, which is an intrinsic geometric flow on a Riemannian manifold, is defined by
\[\frac{\partial}{\partial t}g=-Rg,\]
while the volume preserving normalized Yamabe flow is
\begin{equation}\label{1.2}
\frac{\partial}{\partial t}g=(r-R)g,
\end{equation}
where $r$ is the average scalar curvature. Given an initial metric on a compact locally conformally flat manifold with positive Ricci curvature, Ye \cite{Ye94} proved that the solution $g(t)$ to the normalized Yamabe flow \eqref{1.2} exists for all time and converges in $C^\infty$ norm to a conformal metric of constant scalar curvature.  In 2005, Brendle \cite{Br05} proved the flow \eqref{1.2} on a closed manifold exists for all time and converges to a metric with constant scalar curvature if the dimension $n$ of the underlying manifold is  $3\leq n\leq 5$ or the initial metric is locally conformally flat.

Catino, Cremaschi, Djadli, Mantegazza and Mazzieri \cite{RBflow} proved short time existence and curvature estimates for the Ricci-Bourguignon flow
\[\frac{\partial}{\partial t}g=-2Ric+2\rho Rg,\]
which is a family of geometric flows introduced by Bourguignon \cite{B81} in 1981 and can be viewed as a special case of \eqref{1.1} by taking $\alpha=2$ and $\beta=-2\rho$.

Liang and Zhu \cite{LZ19} proved that the norm of the Weyl tensor of any smooth solution to the Ricci-Bourguignon flow can be explicitly estimated in terms of its initial value on a given ball and a local uniform bound on the Ricci tensor. As an application, they \cite{LZ19} concluded that the Riemannian curvature is uniformly bounded along the Ricci-Bourguignon flow defined on a finite time interval with uniformly bounded Ricci tensor. More recently, Qiu and Zhu \cite{QZ21} showed the short time existence of the Ricci-Bourguignon
flow on a compact Riemannian manifold with constant mean curvature on the boundary if the initial metric has constant mean curvature and satisfies
some compatibility conditions.

In this paper, we consider short and long time existence of the $(\alpha,\beta)$-Ricci-Yamabe flow \eqref{1.1} on a smooth closed Riemannian manifold. The first main result is a short time existence theorem of a smooth solution $g(x,t)$ to \eqref{1.1}.
\begin{theorem}\label{thm1.1}
Let $(M,g_0)$ be an $n$-dimensional $(n\geq2)$ smooth closed Riemannian manifold. Then the $(\alpha,\beta)$-Ricci-Yamabe flow \eqref{1.1} with $\alpha>0$, $\beta>-\frac{\alpha}{n-1}$ and the initial metric $g(\cdot,0)=g_0$ has a smooth solution $g(x,t)$ on $M\times[0,\epsilon)$ for some $\epsilon>0$.
\end{theorem}

There exists a long time existence criteria for smooth solutions to \eqref{1.1} that is a generalization of Hamilton's long time existence theorem for the Ricci flow (see \cite{Ha95}).
\begin{theorem}\label{cor1.1}
For any smooth initial metric $g(\cdot,0)=g_0$ on an $n$-dimensional $(n\geq2)$ closed Riemannian manifold, there exists a maximal time $T$ on which there is a smooth solution $g(x,t)$ for $0\leq t<T$ to the $(\alpha,\beta)$-Ricci-Yamabe flow \eqref{1.1} with $\alpha>0$ and $\beta>-\frac{\alpha}{n-1}$. Either $T=\infty$ or the Riemannian curvature is unbounded as $t\rightarrow T$.
\end{theorem}

Similary, we also have a smooth extension result for a smooth solution $g(x,t)$ to \eqref{1.1} with conditions on the Ricci curvature and the scalar curvature.
%that is a extension of \v{S}e\v{s}um's long time existence theorem for the Ricci flow (see \cite{Se05}).
\begin{theorem}\label{cor1.2}
Let $g(x,t)$, $0\leq t<T<+\infty$, be a smooth solution to the $(\alpha,\beta)$-Ricci-Yamabe flow \eqref{1.1} with $0<\alpha\leq1$ and $\beta>-\frac{\alpha}{n-1}$ on an $n$-dimensional $(n\geq2)$ closed  Riemannian manifold $M$. If
\[\sup_{M\times[0,T)}(|Ric|+|\nabla^2R|)<\infty,\]
then this flow can be extended smoothly over time $T$.
\end{theorem}

Moreover, we derive a long time existence theorem of the  $(\alpha,\beta)$-Ricci-Yamabe flow \eqref{1.1} in terms of its scalar curvature and Weyl tensor, which extends a previous  result of Cao \cite{cao11}.

\begin{theorem}\label{thm1.4}
Let $g(x,t)$, $0\leq t<T<+\infty$, be a maximal smooth solution to the $(\alpha,\beta)$-Ricci-Yamabe flow \eqref{1.1} with $\alpha>0$ and $\beta>-\frac{\alpha}{n-1}$ on an $n$-dimensional $(n\geq3)$ closed  Riemannian manifold $M$. Then either
\[\limsup_{t\rightarrow T}(\max_MR)=\infty\]
or
\[\limsup_{t\rightarrow T}(\max_MR)<\infty\ \text{but}\ \limsup_{t\rightarrow T}\big(\max_M\frac{|W|+|\nabla R|+|\nabla^2 R|}{R}\big)=\infty.\]
\end{theorem}

This paper is arranged as follows. In Section \ref{sec:2}, we establish the short time existence of a smooth solution $g(x,t)$ to the $(\alpha,\beta)$-Ricci-Yamabe flow \eqref{1.1}. We mainly derive evolution equations of curvature operators in Section \ref{sec:3}. In Section \ref{sec:4}, we estimate every ordered covariant derivative of the Riemannian curvature step by step and give a proof of Theorem \ref{cor1.1}. In Section \ref{sec:5}, we show how the Riemannian curvature of $g(x,t)$ can be locally controlled by the Riemannian curvature of the initial metric, the Ricci curvature and the second-order derivative of the scalar curvature. Moreover, we finish the proof of Theorem \ref{cor1.2}. In the last section, we prove Theorem \ref{thm1.4} that relates the long time existence of the $(\alpha,\beta)$-Ricci-Yamabe flow to the behavior of the scalar curvature and Weyl tensor.
\section{Short time existence} % (fold)
\label{sec:2}
In this section, we prove the short time existence of the $(\alpha,\beta)$-Ricci-Yamabe flow \eqref{1.1} by DeTurck's trick which DeTurk \cite{DT83,DT03} established to show that the Ricci flow is equivalent to an initial value problem for a strictly parabolic linear second ordered partial differential equation.
\\\\\textbf{Proof of Theorem \ref{thm1.1}:} For each $h\in\Gamma(S^2M)$, the differentials of the Ricci curvature and the scalar curvature at normal coordinates of $g_0$ in the direction of $h$ (see Theorem 1.174 in \cite{Besse}, \cite{RBflow} or \cite{T06}) are
\[DRic_{g_0}(h)_{ik}=-\frac{1}{2}(\Delta_{g_0}h_{ik}+\nabla_i\nabla_ktr_{g_0}h-\nabla_i\nabla_jh_{kj}-\nabla_k\nabla_jh_{ij})+\cdot\cdot\cdot,\]
and
\[DR_{g_0}(h)=-\Delta_{g_0}tr_{g_0}h+\nabla_j\nabla_lh_{jl}+\cdot\cdot\cdot,\]
respectively, where we omit the lower ordered terms. Then the linearization (without lower ordered terms) of the second ordered nonlinear partial differential operator $P=-2\alpha Ric-\beta Rg$ at $g_0$ is
\begin{eqnarray}\label{3.1}
DP_{g_0}(h)_{ik}&=&-2\alpha DRic_{g_0}(h)_{ik}-\beta DR_{g_0}(h)(g_0)_{ik}-\beta R_{g_0}h_{ik}\notag\\
&=&\alpha(\Delta_{g_0}h_{ik}+\nabla_i\nabla_ktr_{g_0}h-\nabla_i\nabla_jh_{kj}-\nabla_k\nabla_jh_{ij})\notag\\
&&+\beta(\Delta_{g_0}tr_{g_0}h-\nabla_j\nabla_lh_{jl})(g_0)_{ik}+\cdot\cdot\cdot
\end{eqnarray}

Let $X:\Gamma(S^2M)\rightarrow\Gamma(TM)$ be the vector field defined by
%&:=&-\alpha g_0^{jk})h^{pq}\nabla_p(\frac{1}{2}tr_g(g_0)g_{qk}-(g_0)_{qk})\notag\\
\[X^j(g):=-\frac{\alpha}{2}g_0^{jk}g^{pq}(\nabla_k(g_0)_{pq}-\nabla_p(g_0)_{qk}-\nabla_q(g_0)_{pk}.\]
Therefore, the linearization (without lower ordered terms) of the Lie derivative $L_Xg_0$ of $g_0$ in the direction of $X$  is
\begin{equation}\label{3.4}
(DL_X)_{g_0}(h)_{ik}=\alpha(\nabla_i\nabla_ktr_{g_0}h-\nabla_i\nabla_jh_{kj}-\nabla_k\nabla_jh_{ij})+\cdot\cdot\cdot.
\end{equation}

Following  DeTuck's trick (\cite{DT83,DT03}), we only need to check the operator $D(P-L_X)_{g_0}$ is strongly elliptic.

For an arbitrary cotangent vector $\xi$, the principal symbol of the linear differential opeartor $D(P-L_X)_{g_0}(h)_{ik}$ in the direction of $\xi$ is
\begin{eqnarray}\label{3.2}
\sigma_\xi(D(P-L_X)_{g_0})(h)_{ik}&=&\alpha(\xi_j\xi_jh_{ik}+\xi_i\xi_ktr_{g_0}h-\xi_i\xi_jh_{kj}-\xi_k\xi_jh_{ij})\notag\\
&&+\beta(\xi_j\xi_jtr_{g_0}h-\xi_j\xi_lh_{jl})(g_0)_{ik}\notag\\
&&-\alpha(\xi_i\xi_ktr_{g_0}h-\xi_i\xi_jh_{kj}-\xi_k\xi_jh_{ij})\notag\\
&=&\alpha\xi_j\xi_jh_{ik}+\beta(\xi_j\xi_jtr_{g_0}h-\xi_j\xi_lh_{jl})(g_0)_{ik}
\end{eqnarray}

Since \eqref{3.2} is homogeneous, we can choose an orthonormal basis $\{e_i\}_{i=1,...,n}$ at such a point so that $(g_0)_{ik}=\delta_{ik}$ and $|\xi|=1$ with $\xi_1=1$ while $\xi_a=0$ for $a\in\{2,...,n\}$ without loss of generality. Hence,
\begin{equation}\label{3.3}
\sigma_\xi(D(P-L_X)_{g_0})(h)_{ik}=\alpha h_{ik}+\beta(\delta_{ik}tr_{g_0}h-h_{11}\delta_{ik}).
\end{equation}

Similar as in \cite{RBflow}, $\sigma_\xi(D(P-L_X)_{g_0})$ can be represented in the coordinate system
\[(h_{11},h_{22},...,h_{nn},h_{12},...,h_{1n},h_{24},...,h_{n-1,n})\]
by the matrix of
\[
\begin{pmatrix}
U_{n\times n}& \\
 &  \alpha Id_{\frac{n(n-1)}{2}\times\frac{n(n-1)}{2}}
\end{pmatrix},
\]
where all the omitted terms are zero and
\[U_{n\times n}=
\begin{pmatrix}
\alpha & \beta & \beta & \cdots & \beta\\
0 & \alpha+\beta & \beta & \cdots & \beta\\
0 & \beta & \alpha+\beta & \cdots & \beta\\
\vdots &\vdots &\vdots &\ddots &\vdots\\
0 &\beta &\beta &\cdots & \alpha+\beta
\end{pmatrix}.
\]
We denote the submatrix
\[
\begin{pmatrix}
 \alpha+\beta & \beta & \cdots & \beta\\
 \beta & \alpha+\beta & \cdots & \beta\\
\vdots &\vdots &\ddots &\vdots\\
\beta &\beta &\cdots & \alpha+\beta
\end{pmatrix}
\]
in $U_{n\times n}$ by $V_{(n-1)\times (n-1)}$. It is easy to verify that
\[det(V_{(n-1)\times (n-1)}-\lambda Id_{(n-1)\times (n-1)})=(\alpha-\lambda)^{n-2}(\alpha+(n-1)\beta-\lambda)\]
by induction.

Therefore, the eigenvalues of $\sigma_\xi(D(P-L_X)_{g_0})$ are $\alpha>0$ and $\alpha+(n-1)\beta>0$ with multiplicities $\frac{n(n+1)}{2}-1$ and $1$, respectively. This implies that the operator $D(P-L_X)_{g_0}$ is strongly elliptic.

We conclude that there exists a smooth solution $g(x,t)$ to the $(\alpha,\beta)$-Ricci-Yamabe flow \eqref{1.1} with initial metric $g(\cdot,0)=g_0$ on $M\times[0,\epsilon)$ for some $\epsilon>0$.$\hfill\Box$

\section{Curvature evolution equations} % (fold)
\label{sec:3}
In this section, we calculate evolution equations for the Riemannian curvature $Rm$, Ricci curvature $Ric$ and the scalar curvature $R$ of a smooth solution $g(x,t)$ to the $(\alpha,\beta)$-Ricci-Yamabe flow \eqref{1.1} as Hamilton \cite{3D} did for the Ricci flow in 1980s. Moreover, we show some curvature conditions that are preserved along the $(\alpha,\beta)$-Ricci-Yamabe flow \eqref{1.1}.

The notations we used are in accordance with that in Hamilton \cite{3D}. The Riemannian metric is $g_{ij}$ (without the subscripts $x$ and $t$ in the components) and its inverse is $g^{ij}$. The Levi-Civita connection is given by the Christoffel symbols
\[\Gamma_{ij}^h=\frac{1}{2}g^{hk}\big(\frac{\partial}{\partial x^i}g_{jk}+\frac{\partial}{\partial x^j}g_{ik}-\frac{\partial}{\partial x^k}g_{ij}\big).\]

The Riemannian curvature is
\[R_{ijk}^h=\frac{\partial}{\partial x^i}\Gamma_{jk}^h-\frac{\partial}{\partial x^j}\Gamma_{ik}^h+\Gamma_{ip}^h\Gamma_{jk}^p-\Gamma_{jp}^h\Gamma_{ik}^p\]
and
\[R_{ijkl}=g_{hk}R_{ijl}^h.\]
The Ricci curvature is the contraction
\[R_{ik}=g^{jl}R_{ijkl},\]
while the scalar curvature is
\[R=g^{ij}R_{ij}.\]

The tensor introduced by Hamilton \cite{3D}
\[B_{ijkl}=g^{pr}g^{qs}R_{piqj}R_{rksl}\]
satisfies the symmetries
\[B_{ijkl}=B_{jilk}=B_{klij}.\]

In fact, the evolution equations of the Christoffel symbols and the scalar curvature can be found in \cite{RBflow}. For completeness, we give a detail proof in this section.

The following two formulas, which are
independent of any evolution equation, given in Hamilton \cite{3D} are needed in the proof of curvature evolution equations.
\begin{lemma}[Hamilton \cite{3D}]\label{lem2.1}
For any metric $g_{ij}$ the curvature tensor $R_{ijkl}$ satisfies the identity
\begin{eqnarray}\label{2.0}
&&\nabla_i\nabla_kR_{jl}-\nabla_i\nabla_lR_{jk}-\nabla_j\nabla_kR_{il}+\nabla_j\nabla_lR_{ik}\notag\\
&=&\Delta R_{ijkl}+2(B_{ijkl}-B_{ijlk}-B_{iljk}+B_{ikjl})\notag\\
&&-g^{pq}(R_{pjkl}R_{qi}+R_{ipkl}R_{qj}),
\end{eqnarray}
while the tensor $B_{ijkl}$ satisfies the identity
\begin{equation}\label{2.00}
g^{jl}(B_{ijkl}-2B_{ijlk})=0.
\end{equation}
\end{lemma}

\begin{theorem}\label{thm2.1}
The Riemannian curvature satisfies the evolution equation
\begin{eqnarray}\label{2.1}
\frac{\partial}{\partial t}R_{ijkl}&=&\alpha\Delta R_{ijkl}+2\alpha(B_{ijkl}-B_{ijlk}-B_{iljk}+B_{ikjl})\notag\\
&&-\alpha g^{pq}(R_{pjkl}R_{qi}+R_{ipkl}R_{qj}+R_{ijpl}R_{qk}+R_{ijkp}R_{ql})\notag\\
&&+\frac{\beta}{2}(\nabla_i\nabla_kRg_{jl}-\nabla_i\nabla_lRg_{jk}-\nabla_j\nabla_kRg_{il}+\nabla_j\nabla_lRg_{ik})\notag\\
&&-\beta RR_{ijkl}.
\end{eqnarray}
\end{theorem}

\begin{proof}
We calculate the formulas in normal coordinates at any fixed $(x,t)$. Using \eqref{1.1}, we have
\begin{eqnarray}\label{2.2}
\frac{\partial}{\partial t}\Gamma_{jk}^l&=&\frac{1}{2}(\nabla_j\frac{\partial}{\partial t}g_{kl}+\nabla_k\frac{\partial}{\partial t}g_{jl}-\nabla_l\frac{\partial}{\partial t}g_{jk})\notag\\
&=&-\alpha(\nabla_jR_{kl}+\nabla_kR_{jl}-\nabla_lR_{jk})\notag\\
&&-\frac{\beta}{2}(\nabla_jRg_{kl}+\nabla_kRg_{jl}-\nabla_lRg_{jk}).
\end{eqnarray}

It follows from \eqref{2.2} that
\begin{eqnarray}\label{2.3}
\frac{\partial}{\partial t}R_{ijl}^h&=&\nabla_i\frac{\partial}{\partial t}\Gamma_{jl}^h-\nabla_j\frac{\partial}{\partial t}\Gamma_{il}^h\notag\\
&=&-\alpha\nabla_i(\nabla_jR_{hl}+\nabla_lR_{hj}-\nabla_hR_{jl})\notag\\
&&-\frac{\beta}{2}\nabla_i(\nabla_jRg_{hl}+\nabla_lRg_{hj}-\nabla_hRg_{jl})\notag\\
&&+\alpha\nabla_j(\nabla_iR_{hl}+\nabla_lR_{hi}-\nabla_hR_{il})\notag\\
&&+\frac{\beta}{2}\nabla_j(\nabla_iRg_{hl}+\nabla_lRg_{hi}-\nabla_hRg_{il})\notag\\
&=&-\alpha(R_{ijhp}R_{pl}+R_{ijlp}R_{ph})\notag\\
&&+\alpha(\nabla_i\nabla_hR_{jl}-\nabla_i\nabla_lR_{hj}-\nabla_j\nabla_hR_{il}+\nabla_j\nabla_lR_{hi})\notag\\
&&+\frac{\beta}{2}(\nabla_i\nabla_hRg_{jl}-\nabla_i\nabla_lRg_{hj}-\nabla_j\nabla_hRg_{il}+\nabla_j\nabla_lRg_{hi}).
\end{eqnarray}

%Noticing that
%\[\nabla_j\nabla_iR_{hl}-\nabla_i\nabla_jR_{hl}=R_{jihp}R_{pl}+R_{jilp}R_{hp}=-R_{ijhp}R_{pl}+R_{ijpl}R_{hp},\]
%and then using \eqref{2.0}, we can get
Hence,
\begin{eqnarray}\label{2.5}
\frac{\partial}{\partial t}R_{ijkl}&=&\big(\frac{\partial}{\partial t}g_{hk}\big)R_{ijl}^h+g_{hk}\frac{\partial}{\partial t}R_{ijl}^h\notag\\
&=&-2\alpha R_{pk}R_{ijpl} -\beta RR_{ijkl}-\alpha(R_{ijkp}R_{pl}-R_{ijpl}R_{pk})\notag\\
&&+\alpha(\nabla_i\nabla_kR_{jl}-\nabla_i\nabla_lR_{kj}-\nabla_j\nabla_kR_{il}+\nabla_j\nabla_lR_{ik})\notag\\
&&+\frac{\beta}{2}(\nabla_i\nabla_kRg_{jl}-\nabla_i\nabla_lRg_{jk}-\nabla_j\nabla_kRg_{il}+\nabla_j\nabla_lRg_{ik})\notag\\
&=&\alpha\Delta R_{ijkl}+2\alpha(B_{ijkl}-B_{ijlk}-B_{iljk}+B_{ikjl})\notag\\
&&-\alpha (R_{pjkl}R_{pi}+R_{ipkl}R_{pj}+R_{ijpl}R_{pk}+R_{ijkp}R_{pl})\notag\\
&&+\frac{\beta}{2}(\nabla_i\nabla_kRg_{jl}-\nabla_i\nabla_lRg_{jk}-\nabla_j\nabla_kRg_{il}+\nabla_j\nabla_lRg_{ik})\notag\\
&&-\beta RR_{ijkl},
\end{eqnarray}
where we used \eqref{1.1} and \eqref{2.3} in the second equality, and \eqref{2.0} in the second. This proves \eqref{2.1}.
\end{proof}

Then we derive evolution equations for the Ricci curvature and the scalar curvature, respectively.
\begin{theorem}\label{thm2.2}
The Ricci curvature satisfies the evolution equation
\begin{eqnarray}\label{2.6}
\frac{\partial}{\partial t}R_{ik}&=&\alpha\Delta R_{ik}+2\alpha g^{pr}g^{qs}R_{piqk}R_{rs}-2\alpha g^{pq}R_{pi}R_{qk}\notag\\
&&+\frac{\beta}{2}[(n-2)\nabla_i\nabla_kR+\Delta Rg_{ik}],
\end{eqnarray}
while the scalar curvature satisfies
\begin{equation}\label{2.7}
\frac{\partial}{\partial t}R=[\beta(n-1)+\alpha]\Delta R+2\alpha|Ric|^2+\beta R^2,
\end{equation}
where $|Ric|^2=g^{ik}g^{jl}R_{ij}R_{kl}$.
\end{theorem}

\begin{proof}
It follows from \eqref{1.1} that
\begin{equation}\label{2.4}
\frac{\partial}{\partial t}g^{ij}=-g^{ip}g^{jq}\frac{\partial}{\partial t}g_{pq}=2\alpha g^{ip}g^{jq}R_{pq}+\beta Rg^{ij}.
\end{equation}

By direct computation, we have
\begin{eqnarray*}
\frac{\partial}{\partial t}R_{ik}&=&\frac{\partial}{\partial t}(g^{jl}R_{ijkl})\notag\\
&=&-g^{pj}g^{ql}\big(\frac{\partial}{\partial t}g_{pq}\big)R_{ijkl}+g^{jl}\frac{\partial}{\partial t}R_{ijkl}\notag\\
&=&2\alpha g^{pj}g^{ql}R_{pq}R_{ijkl}+\beta RR_{ik}\notag\\
&&+\alpha\Delta R_{ik}+2\alpha g^{jl}(B_{ijkl}-B_{ijlk}-B_{iljk}+B_{ikjl})\notag\\
&&-2\alpha g^{jl}g^{pq}R_{ipkl}R_{qj}-2\alpha g^{pq}R_{pi}R_{qk}\notag\\
&&+\frac{\beta}{2}(n\nabla_i\nabla_kR-2\nabla_i\nabla_kR+\Delta Rg_{ik})-\beta RR_{ik}\notag\\
&=&\alpha\Delta R_{ik}+2\alpha g^{jl}(B_{ijlk}-B_{iljk}+B_{ikjl})\notag\\
&&-2\alpha g^{pq}R_{pi}R_{qk}+\frac{\beta}{2}[(n-2)\nabla_i\nabla_kR+\Delta Rg_{ik}]\notag\\
&=&\alpha\Delta R_{ik}+2\alpha g^{pr}g^{qs}R_{piqk}R_{rs}-2\alpha g^{pq}R_{pi}R_{qk}\notag\\
&&+\frac{\beta}{2}[(n-2)\nabla_i\nabla_kR+\Delta Rg_{ik}],
\end{eqnarray*}
where we used \eqref{2.1} and \eqref{2.4} in the third equality, \eqref{2.00} in the fourth.

Similarly,
\begin{eqnarray*}
\frac{\partial}{\partial t}R&=&\frac{\partial}{\partial t}(g^{ik}R_{ik})\notag\\
&=&-g^{pi}g^{qk}\big(\frac{\partial}{\partial t}g_{pq}\big)R_{ik}+g^{ik}\frac{\partial}{\partial t}R_{ik}\notag\\
&=&2\alpha g^{pi}g^{qk}R_{pq}R_{ik}+\beta R^2+\alpha\Delta R+2\alpha g^{pr}g^{qs}R_{pq}R_{rs}\notag\\
&&-2\alpha g^{ik}g^{pq}R_{pi}R_{qk}+\frac{\beta}{2}[(n-2)\Delta R+n\Delta R]\notag\\
&=&[\beta(n-1)+\alpha]\Delta R+2\alpha|Ric|^2+\beta R^2,
\end{eqnarray*}
where we used \eqref{2.4} in the second equality and \eqref{2.6} in the third.

This completes the proof.
\end{proof}

As applications of \eqref{2.7}, we show how the lower bound of the scalar curvature can be preserved along the $(\alpha,\beta)$-Ricci-Yamabe flow \eqref{1.1}.
\begin{proposition}\label{prp2.1}
On an $n$-dimensional $(n\geq2)$ closed Riemannian manifold $M$, let $g(x,t)$, $(x,t)\in M\times[0,T)$, be a maximal smooth solution to the $(\alpha,\beta)$-Ricci-Yamabe flow \eqref{1.1} with $\alpha>0$ and $\beta>-\frac{\alpha}{n-1}$. If $R_{g(0)}\geq a$ for some constant $a$, then $R_{g(t)}\geq a$ for any $t\in(0,T)$. In particular, if $a>0$, then $T\leq\frac{n(n-1)}{(n-2)a\alpha}$. Moreover, if $a=0$, then $R_{g(t)}>0$ for any $t\in(0,T)$ or $Ric\equiv0$ along this flow.
%the minimum of the scalar curvature $R$ is nondecreasing along this flow. Moreover,
\end{proposition}

\begin{proof}
It follows from \eqref{2.7} that
\begin{eqnarray}\label{2.11}
\frac{\partial}{\partial t}R&=&[\beta(n-1)+\alpha]\Delta R+2\alpha|Ric|^2+\beta R^2\notag\\
&\geq&[\beta(n-1)+\alpha]\Delta R+(\frac{2\alpha}{n}+\beta)R^2\notag\\
&>&[\beta(n-1)+\alpha]\Delta R+\frac{(n-2)\alpha}{n(n-1)} R^2.
\end{eqnarray}

Since $n\geq2$ and $\beta(n-1)+\alpha>0$, the minimum of $R$ on $M$ satisfies

%$\frac{(n-2)\alpha}{n(n-1)} R^2\geq0$. Hence, it is clear that the minimum of the scalar curvature $R$ is nondecreasing along the flow \eqref{1.1} by applying the maximum principle to \eqref{2.11}.
\begin{equation}\label{2.13}
\frac{d}{dt}R_{\min}\geq\frac{(n-2)\alpha}{n(n-1)} R_{\min}^2\geq0\ \text{a.e.},
\end{equation}
which shows that $$R_{g(\cdot,t)}\geq R_{\min}(t)\geq R_{\min}(0)\geq a$$ for any $t\in[0,T)$. This proves the first result.

If $a>0$, we can rewrite \eqref{2.13} to
\begin{equation}\label{2.12}
\frac{d}{dt}\frac{1}{R_{\min}}\leq-\frac{(n-2)\alpha}{n(n-1)}.
\end{equation}
Integrating \eqref{2.12} from $0$ to $t$ for any $0<t<T$ and using the fact of $R_{\min}(0)\geq a$ gives that
\[R_{\min}(t)\geq\frac{n(n-1)a}{n(n-1)-(n-2)a\alpha t},\]
Therefore, the arbitrary of $t$ implies the maximal existence time $T$ satisfies $$T\leq\frac{n(n-1)}{(n-2)a\alpha}.$$

Now we consider the case of $R_{g(\cdot,0)}\geq0$. by applying the strong maximum principle to \eqref{2.11}, we know that $R_{g(\cdot,t)}>0$ or $R_{g(\cdot,t)}\equiv0$ for any $0<t<T$. In the latter case, \eqref{2.7} implies that $Ric_{g(\cdot,t)}\equiv0$.

This completes the proof.
\end{proof}

There is a result for the stationary solution to \eqref{1.1} when it exists.
\begin{theorem}
The stationary solution to the $(\alpha,\beta)$-Ricci-Yamabe flow \eqref{1.1} on an $n$-dimensional $(n\geq3)$ closed Riemannian manifold $M$ with $\alpha>0$ and $\beta>-\frac{\alpha}{n-1}$ must be Ricci-flat.
\end{theorem}

\begin{proof}
Let $g_\infty$ be the stationary solution to \eqref{1.1}. Then its Ricci curvature and scalar curvature satisfy
\begin{equation}\label{2.14}
Ric_{\infty}=-\frac{\beta}{2\alpha}R_{\infty}g_\infty.
\end{equation}

By Schur's lemma (see e.g. \cite{peterson}), $(M,g_\infty)$ must be Einstein. Assume that $R_\infty\neq0$, then
\[-\frac{\beta}{2\alpha}=\frac{1}{n},\]
i.e.,
\[\beta=-\frac{2\alpha}{n}<-\frac{\alpha}{n-1}\]
since $\alpha>0$ and $n\geq3$, which is a contradiction.

Therefore, We conclude $(M,g_\infty)$ is Ricci flat.
\end{proof}
\section{Derivative estimates for the Riemannian curvarure and long time existence \Rmnum{1}} % (fold)
\label{sec:4}
In this section, we show that every ordered covariant derivative of the Riamnnian curvature of smooth solutions to \eqref{1.1} is bounded by only assuming the Riemannian curvature is uniformly bounded. As an application, we prove a long time existence theorem (Theorem \ref{cor1.1}) for \eqref{1.1}.

As in \cite{Shi89,LYZ20,KMW16}, we use the standard $*$-notation that $A*B$ represents some linear combination of contractions of the tensor product $A\otimes B$ of any time-dependent tensor fields $A$ and $B$ by using the metric $g(x,t)$.

Hence, \eqref{2.1} and \eqref{2.2} can be rewritten as
\begin{equation}\label{4.1}
\frac{\partial}{\partial t}Rm=\alpha \Delta Rm-\beta R\cdot Rm+\alpha g^{-1}*g^{-1}*Rm*Rm+\beta\nabla^2R*g,
\end{equation}
and
%\begin{equation}\label{4.2}
\[\frac{\partial}{\partial t}\Gamma=g^{-1}*\nabla(\alpha g^{-1}*Rm+\beta Rg),\]
%\end{equation}
respectively.

More generally, a time-dependent tensor field $A$ satisfies
\begin{eqnarray}\label{4.3}
\frac{\partial}{\partial t}\nabla A&=&\nabla(\frac{\partial}{\partial t}A)+(\frac{\partial}{\partial t}\Gamma)*A\notag\\
&=&\nabla(\frac{\partial}{\partial t}A)+g^{-1}*\nabla(\alpha g^{-1}*Rm+\beta Rg)*A,
\end{eqnarray}
and
\begin{equation}\label{4.4}
\nabla\Delta A=\Delta\nabla A+g^{-1}*g^{-1}*Rm*\nabla A+g^{-1}*g^{-1}*\nabla Rm*A.
\end{equation}

Then we prove the evolution equation for $k$th-order covariant derivative of the Riemannian curvature.
\begin{proposition}\label{prp4.1}
For any nonnegative integer $k$, we have
\begin{eqnarray}\label{4.5}
\frac{\partial}{\partial t}\nabla^kRm&=&\alpha\Delta\nabla^kRm+\alpha\sum_{i+j=k}\nabla^iRm*\nabla^jRm*g^{-1}*g^{-1}\notag\\
&&+\beta\nabla^{k+2}R*g+\beta\sum_{i+j=k}\nabla^iRg*\nabla^jRm*g^{-1}
\end{eqnarray}
and
\begin{eqnarray}\label{4.2}
\frac{\partial}{\partial t}\nabla^kR&=&[\beta(n-1)+\alpha](\Delta\nabla^kR+\sum_{i+j=k}\nabla^iRm*\nabla^jR*g^{-1}*g^{-1})\notag\\
&&+\alpha\sum_{i+j=k}\nabla^iRic*\nabla^jRic+\beta\sum_{i+j=k}\nabla^iR*\nabla^jR
\end{eqnarray}
along the $(\alpha,\beta)$-Ricci-Yamabe flow \eqref{1.1}.
\end{proposition}

\begin{proof}
We prove this result by induction. The case of $k=0$ holds since \eqref{2.1} and \eqref{2.7}, respectively. Suppose that \eqref{4.5} and \eqref{4.2} hold for $k=m>0$, we need to show that they hold for $k=m+1$.

Substituting $A$ by $\nabla^mRm$ in \eqref{4.3} and \eqref{4.4}, and using the inductive hypothesis, we obtain that
\begin{eqnarray*}
\frac{\partial}{\partial t}\nabla^{m+1}Rm&=&\nabla\Big(\alpha\Delta\nabla^mRm+\alpha\sum_{i+j=m}\nabla^iRm*\nabla^jRm*g^{-1}*g^{-1}\notag\\
&&+\beta\nabla^{m+2}R*g+\beta\sum_{i+j=m}\nabla^iRg*\nabla^jRm*g^{-1}\Big)\notag\\
&&+g^{-1}*\nabla(\alpha g^{-1}*Rm+\beta Rg)*\nabla^mRm\notag\\
&=&\alpha\nabla\Delta\nabla^mRm+\alpha\sum_{i+j=m+1}\nabla^iRm*\nabla^jRm*g^{-1}*g^{-1}\notag\\
&&+\beta\nabla^{m+3}R*g+\beta\sum_{i+j=m+1}\nabla^iRg*\nabla^jRm*g^{-1}\notag\\
&=&\alpha\Delta\nabla^{m+1}Rm+\alpha g^{-1}*g^{-1}*Rm*\nabla^{m+1}Rm\notag\\
&&+\alpha g^{-1}*g^{-1}*\nabla Rm*\nabla^mRm+\alpha\sum_{i+j=m+1}\nabla^iRm*\nabla^jRm*g^{-1}*g^{-1}\notag\\
&&+\beta\nabla^{m+3}R*g+\beta\sum_{i+j=m+1}\nabla^iRg*\nabla^jRm*g^{-1}\notag\\
&=&\alpha\Delta\nabla^{m+1}Rm+\alpha\sum_{i+j=m+1}\nabla^iRm*\nabla^jRm*g^{-1}*g^{-1}\notag\\
&&+\beta\nabla^{m+3}R*g+\beta\sum_{i+j=m+1}\nabla^iRg*\nabla^jRm*g^{-1},
\end{eqnarray*}
which yields \eqref{4.5} for $k=m+1$.

Similarly, we have
\begin{eqnarray*}
\frac{\partial}{\partial t}\nabla^{m+1}R&=&\nabla\Big\{[\beta(n-1)+\alpha](\Delta\nabla^mR+\sum_{i+j=m}\nabla^iRm*\nabla^jR*g^{-1}*g^{-1})\notag\\
&&+\alpha\sum_{i+j=m}\nabla^iRic*\nabla^jRic+\beta\sum_{i+j=m}\nabla^iR*\nabla^jR\Big\}\notag\\
&&+g^{-1}*\nabla(\alpha g^{-1}*Rm+\beta Rg)*\nabla^mR\notag\\
&=&[\beta(n-1)+\alpha]\big(\Delta\nabla^{m+1}R+g^{-1}*g^{-1}*Rm*\nabla^{m+1}R\notag\\
&&+g^{-1}*g^{-1}*\nabla Rm*\nabla^mR+\sum_{i+j=m+1}\nabla^iRm*\nabla^jR*g^{-1}*g^{-1}\big)\notag\\
&&+\alpha\sum_{i+j=m+1}\nabla^iRic*\nabla^jRic+\beta\sum_{i+j=m+1}\nabla^iR*\nabla^jR\notag\\
&=&[\beta(n-1)+\alpha](\Delta\nabla^{m+1}R+\sum_{i+j=m+1}\nabla^iRm*\nabla^jR*g^{-1}*g^{-1})\notag\\
&&+\alpha\sum_{i+j=m+1}\nabla^iRic*\nabla^jRic+\beta\sum_{i+j=m+1}\nabla^iR*\nabla^jR
\end{eqnarray*}
which yields \eqref{4.2} for $k=m+1$.
\end{proof}

\begin{lemma}\label{lem4.1}
For any nonnegative integer $k$, we have
\begin{eqnarray}\label{4.6}
&&\frac{\partial}{\partial t}|\nabla^kRm|^2\notag\\
&\leq&\alpha\Delta|\nabla^kRm|^2-2\alpha|\nabla^{k+1}Rm|^2+\alpha\sum_{i+j=k}|\nabla^iRm|\cdot|\nabla^jRm|\cdot|\nabla^kRm|\notag\\
&&+|\beta|C_k|\nabla^{k+2}R|\cdot|\nabla^kRm|+|\beta|C_k\sum_{i+j=k}|\nabla^iR|\cdot|\nabla^jRm|\cdot|\nabla^kRm|
\end{eqnarray}
and
\begin{eqnarray}\label{4.7}
&&\frac{\partial}{\partial t}|\nabla^kR|^2\notag\\
&\leq&[\beta(n-1)+\alpha]\Delta|\nabla^kR|^2-2[\beta(n-1)+\alpha]|\nabla^{k+1}R|^2\notag\\
&&+(\alpha+|\beta|)C_k\sum_{i+j=k}|\nabla^iRm|\cdot|\nabla^jRm|\cdot|\nabla^kR|
\end{eqnarray}
along the $(\alpha,\beta)$-Ricci-Yamabe flow \eqref{1.1} with $\alpha>0$ and $\beta>-\frac{\alpha}{n-1}$, where $C_k$ is a positive constant depends only on $n$ and $k$.
\end{lemma}

\begin{proof}
Calculating directly, we have
\begin{eqnarray*}
&&\frac{\partial}{\partial t}|\nabla^kRm|^2\notag\\
&\leq&C_k|\alpha Ric+\beta Rg|\cdot|\nabla^kRm|^2+2\langle\nabla^kRm,\frac{\partial}{\partial t}\nabla^kRm\rangle\notag\\
&\leq&\alpha C_k|Ric|\cdot|\nabla^kRm|^2+|\beta|C_k|R|\cdot|\nabla^kRm|^2+2\alpha\langle\nabla^kRm,\Delta\nabla^kRm\rangle\notag\\
&&+\alpha C_k\sum_{i+j=k}|\nabla^iRm|\cdot|\nabla^jRm|\cdot|\nabla^kRm|\notag\\
&&+|\beta|C_k|\nabla^{k+2}R|\cdot|\nabla^kRm|+|\beta|C_k\sum_{i+j=k}|\nabla^iR|\cdot|\nabla^jRm|\cdot|\nabla^kRm|\notag\\
&\leq&\alpha\Delta|\nabla^kRm|^2-2\alpha|\nabla^{k+1}Rm|^2+\alpha C_k\sum_{i+j=k}|\nabla^iRm|\cdot|\nabla^jRm|\cdot|\nabla^kRm|\notag\\
&&+|\beta|C_k|\nabla^{k+2}R|\cdot|\nabla^kRm|+|\beta|C_k\sum_{i+j=k}|\nabla^iR|\cdot|\nabla^jRm|\cdot|\nabla^kRm|,
\end{eqnarray*}
where we used \eqref{1.1} in the first inequality and \eqref{4.5} in the second. This proves \eqref{4.6}.

Similarly, by using \eqref{1.1}, \eqref{4.2} and Cauchy's inequality, we can obtain that
\begin{eqnarray*}
\frac{\partial}{\partial t}|\nabla^kR|^2&\leq&C_k|\alpha Ric+\beta Rg|\cdot|\nabla^kR|^2+2\langle\nabla^kR,\frac{\partial}{\partial t}\nabla^kR\rangle\notag\\
&\leq&\alpha C_k|Ric|\cdot|\nabla^kR|^2+|\beta|C_k|R|\cdot|\nabla^kR|^2\notag\\
&&+2[\beta(n-1)+\alpha]\langle\nabla^kR,\Delta\nabla^k R\rangle\notag\\
&&+(\alpha+|\beta|)C_k\sum_{i+j=k}|\nabla^iRm|\cdot|\nabla^jR|\cdot|\nabla^kR|\notag\\
&&+\alpha C_k\sum_{i+j=k}|\nabla^iRic|\cdot|\nabla^jRic|\cdot|\nabla^kR|\notag\\
&&+|\beta|C_k\sum_{i+j=k}|\nabla^iR|\cdot|\nabla^jR|\cdot|\nabla^kR|\notag\\
&\leq&[\beta(n-1)+\alpha]\Delta|\nabla^kR|^2-2[\beta(n-1)+\alpha]|\nabla^{k+1}R|^2\notag\\
&&+(\alpha+|\beta|)C_k\sum_{i+j=k}|\nabla^iRm|\cdot|\nabla^jRm|\cdot|\nabla^kR|,
\end{eqnarray*}
which yields \eqref{4.7}.
\end{proof}

Then we show that all the $L^2$-norms of derivatives of the Riemannian curvature and the scalar curvature are bounded if the Riemannian curvature is uniformly bounded. The arguments are motivated by curvature estimates for the Ricci flow by Shi \cite{Shi89} and a geometric flow over K\"{a}hler manifold by Li, Yuan and Zhang \cite{LYZ20}.

The following analytic result is necessary.
\begin{proposition}[Catino et al. \cite{RBflow}]\label{prp4.2}
Let $k\in\mathbb{N}$, $p\in[1,+\infty]$ and $q\in[1,+\infty)$. There exists a constant $C(n,k,p,q)$ such that for all $0\leq j\leq k$ and all tensors $T$
\[\|\nabla^jT\|_{r_j}\leq C\|T\|_p^{1-j/k}\|\nabla^kT\|_q^{j/k},\]
where $1/r_j=(1-j/k)/p+j/k/q$.
\end{proposition}

\begin{theorem}\label{thm4.1}
Let $g(t)$ be a smooth solution to the $(\alpha,\beta)$-Ricci-Yamabe flow \eqref{1.1} with $\alpha>0$ and $\beta>-\frac{\alpha}{n-1}$ on an $n$-dimensional closed Riemannian manifold.
If
\[|Rm|\leq K\]
on $M\times[0,a/K]$ for some positive constants $a$ and $K$, then for each positive integer $k$, there exists a nonnegative constant $B_k$ that depends only on $k$, $n$, $\alpha$, $\beta$, $a$, $K$ and $Vol_{g(0)}(M)$ so that
\begin{equation}\label{4.8}
\|\nabla^kRm\|_{2}^2+\|\nabla^{k+1}R\|_{2}^2\leq\frac{B_k}{t^k}
\end{equation}
on $M$ for all $t\in[0,a/K]$. Here
\[\|A\|_{p}:=\big(\int_M|A|^2dVol_{g(t)}\big)^\frac{1}{p}\]
for any time-dependent $L^p$ function $A$ on $M$ with $1\leq p<\infty$.
\end{theorem}

\begin{proof}
We derive \eqref{4.8} step by step. First of all, we consider the case of $k=1$. Define
\[u_0=t(\|Rm\|_2^2+\|\nabla R\|_2^2)+A_0\|R\|_2^2\]
with $A_0$ is a positive constant to be determined later. We need to compute $\frac{du_0}{dt}$.

It follows from \eqref{1.1} that
\begin{eqnarray}\label{4.9}
\frac{d}{d t}d Vol_{g(t)}&=&\frac{1}{2}g^{ij}\frac{\partial}{\partial t}g_{ij}dVol_{g(t)}\notag\\
&=&-(\alpha+\frac{n\beta}{2})RdVol_{g(t)}\notag\\
&\leq&(n\alpha+\frac{n^2}{2}|\beta|)KdVol_{g(t)},
\end{eqnarray}
i.e.,
\begin{equation}\label{4.10}
\frac{d}{d t}\ln d  Vol_{g(t)}\leq (n\alpha+\frac{n^2}{2}|\beta|)K.
\end{equation}
Integrating \eqref{4.10} from $0$ to $t$, we arrive at
\begin{equation}\label{4.11}
Vol_{g(t)}(M)\leq Vol_{g(0)}(M)\cdot\exp\{(n\alpha+\frac{n^2}{2}|\beta|)Kt\}.
\end{equation}

Lemma \ref{lem4.1} implies that
\begin{equation}\label{4.12}
\frac{\partial}{\partial t}|Rm|^2\leq\alpha\Delta|Rm|^2-2\alpha|\nabla Rm|^2+(\alpha+|\beta|) C_0K^3+|\beta|C_0K|\nabla^{2}R|,
\end{equation}
and
\begin{eqnarray}\label{4.13}
&&\frac{\partial}{\partial t}|\nabla R|^2\notag\\
&\leq&[\beta(n-1)+\alpha]\Delta|\nabla R|^2-2[\beta(n-1)+\alpha]|\nabla^{2}R|^2\notag\\
&&+(\alpha+|\beta|) C_1K|\nabla Rm|\cdot|\nabla R|.
\end{eqnarray}

%Integrating \eqref{4.12} and \eqref{4.13} on $M$, respectively, and using Stokes's formula, we have
Hence,
\begin{eqnarray}\label{4.14}
&&\frac{d}{d t}\int_Mt|Rm|^2dVol_{g(t)}\notag\\
&=&\int_M|Rm|^2dVol_{g(t)}+\int_Mt\frac{\partial}{\partial t}|Rm|^2dVol_{g(t)}+\int_Mt|Rm|^2\frac{d}{d t}dVol_{g(t)}\notag\\
&\leq&K^2Vol_{g(t)}(M)-2\alpha t\int_M |\nabla Rm|^2dVol_{g(t)}+(\alpha+|\beta|)C_0K^3tVol_{g(t)}(M)\notag\\
&&+|\beta|C_0Kt\int_M|\nabla^2R|dVol_{g(t)},
\end{eqnarray}
where we used \eqref{4.9}, \eqref{4.12} and Stokes's formula, and
\begin{eqnarray}\label{4.15}
&&\frac{d}{d t}\int_Mt|\nabla R|^2dVol_{g(t)}\notag\\
&=&\int_M|\nabla R|^2dVol_{g(t)}+\int_Mt\frac{\partial}{\partial t}|\nabla R|^2dVol_{g(t)}+\int_Mt|\nabla R|^2\frac{d}{d t}dVol_{g(t)}\notag\\
&\leq&\int_M|\nabla R|^2dVol_{g(t)}-2[\beta(n-1)+\alpha] t\int_M |\nabla^2 R|^2dVol_{g(t)}\notag\\
&&+(\alpha+|\beta|)C_1Kt\int_M |\nabla Rm|\cdot|\nabla R|dVol_{g(t)},
\end{eqnarray}
where we used \eqref{4.9} and \eqref{4.13} and Stokes's formula.

Moreover, it follows from \eqref{2.7} that
\begin{eqnarray*}
\frac{\partial}{\partial t}R^2&=&[\beta(n-1)+\alpha](\Delta R^2-2|\nabla R|^2)+4\alpha R|Ric|^2+2\beta R^3\notag\\
&\leq&[\beta(n-1)+\alpha](\Delta R^2-2|\nabla R|^2)+(\alpha+|\beta|)C_0K^3.
\end{eqnarray*}
Combining with \eqref{4.9} and using Stokes's formula, we have
\begin{eqnarray}\label{4.16}
\frac{\partial}{\partial t}\int_M A_0R^2dVol_{g(t)}&\leq&-2[\beta(n-1)+\alpha]A_0\int_M |\nabla R|^2dVol_{g(t)}\notag\\
&&+(\alpha+|\beta|)C_0K^3A_0Vol_{g(t)}(M).
\end{eqnarray}

Therefore, \eqref{4.14}, \eqref{4.15} together with \eqref{4.16} gives
\begin{eqnarray}\label{4.17}
\frac{d}{dt}u_0&\leq&-2\alpha t\int_M |\nabla Rm|^2dVol_{g(t)}+(\alpha+|\beta|) C_1Kt\int_M |\nabla Rm|\cdot|\nabla R|dVol_{g(t)}\notag\\
&&+|\beta|C_0Kt\int_M |\nabla^2 R|dVol_{g(t)}-2[\beta(n-1)+\alpha] t\int_M |\nabla^2 R|^2dVol_{g(t)}\notag\\
&&+[1-2(\beta(n-1)+\alpha)A_0]\int_M |\nabla R|^2dVol_{g(t)}\notag\\
&&+[(\alpha+|\beta|)C_0K^3(t+A_0)+K^2]Vol_{g(t)}(M).
\end{eqnarray}

By Cauchy's inequality and the fact of $t\leq a/K$, we obtain that
\begin{eqnarray}\label{4.18}
&&-2\alpha t\int_M |\nabla Rm|^2dVol_{g(t)}+(\alpha+|\beta|) C_1Kt\int_M |\nabla Rm|\cdot|\nabla R|dVol_{g(t)}\notag\\
&\leq& \frac{(\alpha+|\beta|)^2}{\alpha}C_1K^2 t\int_M |\nabla R|^2dVol_{g(t)}\notag\\
&\leq& \frac{(\alpha+|\beta|)^2}{\alpha}C_1 aK\int_M |\nabla R|^2dVol_{g(t)},
\end{eqnarray}
and
\begin{eqnarray}\label{4.19}
&&|\beta|C_0Kt\int_M |\nabla^2 R|dVol_{g(t)}-2[\beta(n-1)+\alpha] t\int_M |\nabla^2 R|^2dVol_{g(t)}\notag\\
&\leq&\frac{\beta^2C_0K^2t}{\beta(n-1)+\alpha}Vol_{g(t)}(M)\leq\frac{\beta^2C_0aK}{\beta(n-1)+\alpha}Vol_{g(t)}(M).
\end{eqnarray}

Plugging \eqref{4.18} and \eqref{4.19} into \eqref{4.17} and then using $t\leq a/K$ and \eqref{4.11}, we get
\begin{eqnarray}\label{4.20}
&&\frac{d}{dt}u_0\notag\\
&\leq&[1+\frac{(\alpha+|\beta|)^2}{\alpha}C_1 aK-2(\beta(n-1)+\alpha)A_0]\int_M |\nabla R|^2dVol_{g(t)}\notag\\
&&+\{(\alpha+|\beta|)C_0(aK^2+A_0K^3)+K^2+\frac{\beta^2C_0aK}{\beta(n-1)+\alpha}\}\notag\\
&&\times\exp\{(n\alpha+\frac{n^2}{2}|\beta|)a\} Vol_{g(0)}(M).
\end{eqnarray}

Choose
\begin{equation}\label{4.21}
A_0=[2\beta(n-1)+2\alpha]^{-1}[1+\frac{(\alpha+|\beta|)^2}{\alpha}C_1a+C_1\alpha aK],
\end{equation}
then the first term of RHS of \eqref{4.20} vanishes. For an arbitary $0<t\leq a/K$, integrating \eqref{4.20} from $0$ to $t$, we conclude that $u_0\leq B_0t$ for some constant $B_0$ depends only on $n$, $\alpha$, $\beta$, $a$, $K$ and $Vol_{g(0)}(M)$. Moreover,
\begin{equation}\label{4.22}
\|Rm\|_2^2+\|\nabla R\|_2^2\leq B_0.
\end{equation}

Next, we estimate $\|\nabla Rm\|_2^2+\|\nabla^2 R\|_2^2$ based on \eqref{4.22}. Define
\[u_1=t^2(\|\nabla Rm\|_2^2+\|\nabla^2 R\|_2^2)+A_1t(\|Rm\|_2^2+\|\nabla R\|_2^2)+A_0\|R\|_2^2,\]
where $A_0$ is defined in \eqref{4.21} and $A_1$ is a positive constant to be determined later.

Then we calculate $\frac{d u_1}{dt}$. It follows from Lemma \ref{lem4.1}, \eqref{4.9} and Stokes's formula that
\begin{eqnarray}\label{4.23}
&&\frac{\partial}{\partial t}\int_M t^2|\nabla Rm|^2dVol_{g(t)}\notag\\
&\leq&-2\alpha t^2\int_M |\nabla^2Rm|^2dVol_{g(t)}+(\alpha+|\beta|)C_1Kt^2\int_M |\nabla Rm|^2dVol_{g(t)}\notag\\
&&+|\beta|C_1t^2\int_M |\nabla^3R|\cdot|\nabla Rm|dVol_{g(t)}+2t\int_M |\nabla Rm|^2dVol_{g(t)},
\end{eqnarray}
and
\begin{eqnarray}\label{4.24}
&&\frac{\partial}{\partial t}\int_M t^2|\nabla^2 R|^2dVol_{g(t)}\notag\\
&\leq&-2[\beta(n-1)+\alpha]t^2\int_M |\nabla^3R|^2dVol_{g(t)}+(\alpha+|\beta|)C_2Kt^2\int_M |\nabla^2 Rm|\cdot|\nabla^2R|dVol_{g(t)}\notag\\
&&+(\alpha+|\beta|)C_2t^2\int_M |\nabla Rm|^2\cdot|\nabla^2R|dVol_{g(t)}+2t\int_M |\nabla^2 R|^2dVol_{g(t)}.
\end{eqnarray}

Combining \eqref{4.14}, \eqref{4.15}, \eqref{4.16}, \eqref{4.23} with \eqref{4.24}, we arrive at
\begin{eqnarray}\label{4.25}
&&\frac{d}{dt}u_1\notag\\
&\leq&-2\alpha t^2\int_M |\nabla^2Rm|^2dVol_{g(t)}+(\alpha+|\beta|)C_2Kt^2\int_M |\nabla^2Rm|\cdot|\nabla^2R|dVol_{g(t)}\notag\\
&&+|\beta|C_1t^2\int_M |\nabla^3R|\cdot|\nabla Rm|dVol_{g(t)}-2[\beta(n-1)+\alpha]t^2\int_M |\nabla^3R|^2dVol_{g(t)}\notag\\
&&+(\alpha+|\beta|)C_1t^2\int_M |\nabla Rm|^2|\nabla^2R|dVol_{g(t)}+A_1(\alpha+|\beta|) C_1Kt\int_M |\nabla Rm|\cdot|\nabla R|dVol_{g(t)}\notag\\
&&+A_1|\beta|C_0Kt\int_M |\nabla^2R|dVol_{g(t)}+[(\alpha+|\beta|)C_1Kt^2+2t-2A_1\alpha t]\notag\\
&&\times\int_M|\nabla Rm|^2dVol_{g(t)}+2[1-A_1(\beta(n-1)+\alpha)]t\int_M |\nabla^2R|^2dVol_{g(t)}\notag\\
&&+[A_1(\alpha+|\beta|)C_1Kt+A_1-2(\beta(n-1)+\alpha)A_0]\int_M |\nabla R|^2Vol_{g(t)}\notag\\
&&+[(\alpha+|\beta|)C_0K^3(A_1t+A_0)+A_1K^2] Vol_{g(t)}(M).
\end{eqnarray}

In the following, we deal with RHS of \eqref{4.25}.

By Cauchy's inequality, we have
\begin{eqnarray}\label{4.26}
&&-2\alpha t^2\int_M |\nabla^2Rm|^2dVol_{g(t)}+(\alpha+|\beta|)C_2Kt^2\int_M |\nabla^2Rm|\cdot|\nabla^2R|dVol_{g(t)}\notag\\
&\leq&\frac{(\alpha+|\beta|)^2C_2K^2t^2}{\alpha}\int_M |\nabla^2R|^2dVol_{g(t)},
\end{eqnarray}
\begin{eqnarray}\label{4.27}
&&|\beta|C_1t^2\int_M |\nabla^3R|\cdot|\nabla Rm|dVol_{g(t)}-2[\beta(n-1)+\alpha]t^2\int_M |\nabla^3R|^2dVol_{g(t)}\notag\\
&\leq&\frac{\beta^2C_1t^2}{\beta(n-1)+\alpha}\int_M |\nabla Rm|^2dVol_{g(t)},
\end{eqnarray}
\begin{eqnarray}\label{4.28}
&&A_1(\alpha+|\beta|)C_1Kt\int_M |\nabla Rm|\cdot|\nabla R|dVol_{g(t)}\notag\\
&\leq&A_1(\alpha+|\beta|) C_1t\int_M |\nabla Rm|^2dVol_{g(t)}\notag\\
&&+A_1(\alpha+|\beta|) C_1K^2t\int_M |\nabla R|^2dVol_{g(t)},
\end{eqnarray}
and
\begin{eqnarray}\label{4.29}
&&A_1|\beta|C_0Kt\int_M |\nabla^2R|dVol_{g(t)}\notag\\
&\leq&A_1[\beta(n-1)+\alpha]t\int_M |\nabla^2 R|^2dVol_{g(t)}+\frac{A_1\beta^2C_0K^2t}{\beta(n-1)+\alpha}Vol_{g(t)}(M).
\end{eqnarray}

Furthermore,
\begin{eqnarray}\label{4.30}
&&(\alpha+|\beta|)C_2t^2\int_M |\nabla Rm|^2|\nabla^2R|dVol_{g(t)}\notag\\
&\leq&(\alpha+|\beta|)C_2t^2\|\nabla Rm\|_4^2\cdot\|\nabla^2R\|_2\notag\\
&\leq&(\alpha+|\beta|)C_2t^2\|Rm\|_\infty\cdot\|\nabla^2Rm\|_2\cdot\|\nabla^2R\|_2\notag\\
&\leq&\alpha t^2\int_M |\nabla^2Rm|^2dVol_{g(t)}+\frac{(\alpha+|\beta|)^2C_2K^2t^2}{\alpha}\int_M |\nabla^2R|^2dVol_{g(t)},
\end{eqnarray}
where we used H\"older's inequality in the first inequality, Proposition \ref{prp4.2} by taking $j=1$, $k=2$, $p=\infty$ and $q=2$ in the second inequality and Cauchy's inequality in the last.

Plugging \eqref{4.26} to \eqref{4.30} into \eqref{4.25} and rearranging, we can obtain
\begin{eqnarray}\label{4.31}
&&\frac{d}{dt}u_1\notag\\
&\leq&\|\nabla Rm\|_2^2\big[(\alpha+|\beta|)C_1Kt+2+\frac{\beta^2C_1t}{\beta(n-1)+\alpha}-A_1(\alpha+|\beta|)\big]t\notag\\
&&+\|\nabla^2R\|_2^2\big[\frac{(\alpha+|\beta|)^2C_2K^2t}{\alpha}+2-A_1(\beta(n-1)+\alpha)\big]t\notag\\
&&+\|\nabla R\|_2^2[A_1(\alpha+|\beta|)C_1K(K+1)t+A_1-2(\beta(n-1)+\alpha)A_0]\notag\\
&&+Vol_{g(t)}(M)\big[(\alpha+|\beta|)C_0K^3(A_1t+A_0)+A_1K^2+\frac{A_1\beta^2C_0K^2t}{\beta(n-1)+\alpha}\big].
\end{eqnarray}

Choose
\[A_1:=\max\{C_1a+\frac{2}{\alpha+|\beta|}+\frac{\beta^2C_1a}{ K[\beta(n-1)+\alpha](\alpha+|\beta|)},\frac{(\alpha+|\beta|)^2C_2Ka}{\alpha[\beta(n-1)+\alpha]}+\frac{2}{\beta(n-1)+\alpha}\}.\]
Applying \eqref{4.22} to \eqref{4.31} and using $t\leq a/K$ and \eqref{4.11} again, we conclude that
\begin{eqnarray}\label{4.32}
&&\frac{d}{dt}u_1\notag\\
&\leq&B_0A_1[(\alpha+|\beta|)C_1(K+1)a+1]+Vol_{g(0)}(M)\exp\{(n\alpha+\frac{n^2}{2}|\beta|)a\}\notag\\
&&\times\big[(\alpha+|\beta|)C_0K^2(A_1a+A_0K)+A_1K^2+\frac{A_1\beta^2C_0Ka}{\beta(n-1)+\alpha}\big].
\end{eqnarray}

For an arbitrary $0<t\leq a/K$, integrating \eqref{4.32} from $0$ to $t$ gives $u_1\leq B_1t$ for some constant $B_1$ depends only on $n$, $\alpha$, $\beta$, $a$, $K$ and $Vol_{g(0)}(M)$.
Moreover,
\begin{equation}\label{4.33}
\|\nabla Rm\|_2^2+\|\nabla^2R\|_2^2\leq\frac{B_1}{t}.
\end{equation}

Define
\[u_2:=t^{3}(\|\nabla^2 Rm\|_2^2+\|\nabla^3R\|_2^2)+\sum_{i=0}^2A_it^i(\|\nabla^{i-1}Rm\|_2^2+\|\nabla^iR\|_2^2).\]
Using \eqref{4.22}, \eqref{4.33} and reasoning similar to \eqref{4.25}, \eqref{4.31}, we can estimate $\frac{d u_2}{dt}$ and then obtain that
\[\|\nabla^2 Rm\|_2^2+\|\nabla^3R\|_2^2\leq\frac{B_2}{t^2}\]
for some constant $B_2$ depends only on $n$, $\alpha$, $\beta$, $a$, $K$ and $Vol_{g(0)}(M)$.

Repeating this procedure step by step, we conclude \eqref{4.8} for each positive integer $k$.
\end{proof}

In the rest of this section, we prove derivative estimates for the Riemannian curvature under the flow \eqref{1.1}.
\begin{theorem}\label{thm4.2}
Let $g(t)$ be a smooth solution to the $(\alpha,\beta)$-Ricci-Yamabe flow \eqref{1.1} with $\alpha>0$ and $\beta>-\frac{\alpha}{n-1}$ on an $n$-dimensional closed Riemannian manifold.
If
\[|Rm|\leq K\]
on $M\times[0,a/K]$ for some positive constants $a$ and $K$, then for each positive integer $k$, there exists a nonnegative constant $\bar{B}_k$ that depends only on $k$, $n$, $\alpha$, $\beta$, $a$, $K$ and $Vol_{g(0)}(M)$ so that
\begin{equation}\label{4.34}
\sup_M|\nabla^kRm|\leq\frac{\bar{B}_k}{t^\frac{k+1}{2}}.
\end{equation}
\end{theorem}

\begin{proof}
For any given $k\geq1$, choose the integer $s=\lceil\frac{n(k+1)}{2}\rceil$ that is the smallest integer larger than $\frac{n(k+1)}{2}$. By taking $p=\infty$ and $q=2$ in Proposition \ref{prp4.2} and then using the fact from Theorem \ref{thm4.1} of
\[\|\nabla^sRm\|_2^2\leq\frac{B_s}{t^s},\]
we have
\begin{equation}\label{4.35}
\|\nabla^{k+1}Rm\|_{\frac{2s}{k+1}}\leq C_kK^{1-\frac{k+1}{s}}\|\nabla^sRm\|_2^\frac{k+1}{s}\leq C_kK^{1-\frac{k+1}{s}}\big(\frac{B_s}{t^s}\big)^\frac{k+1}{2s}.
\end{equation}

By Sobolev embedding theorem (see e.g. Theorem 11.1.1 in \cite{jostpde}), we know that
\[\sup_M|\nabla^kRm|\leq C_kVol_{g(t)}(M)^{\frac{1}{n}-\frac{k+1}{2s}}\|\nabla|\nabla^kRm|\|_\frac{2s}{k+1}.\]

Note that
\[2|\nabla^kRm|\cdot|\nabla|\nabla^kRm||=|\nabla|\nabla^kRm|^2|\leq2|\nabla^kRm|\cdot|\nabla^{k+1}Rm|,\]
i.e.,
\[|\nabla|\nabla^kRm||\leq |\nabla^{k+1}Rm|.\]

Hence,
\begin{eqnarray}\label{4.36}
\sup_M|\nabla^kRm|&\leq&C_kVol_{g(t)}(M)^{\frac{1}{n}-\frac{k+1}{2s}}\|\nabla^{k+1}Rm\|_\frac{2s}{k+1}\notag\\
&\leq&C_kVol_{g(0)}(M)^{\frac{1}{n}-\frac{k+1}{2s}}\cdot\exp\{(\frac{1}{n}-\frac{k+1}{2s})(n\alpha+\frac{n^2}{2}|\beta|)Kt\}\notag\\
&&\times K^{1-\frac{k+1}{s}}\big(\frac{B_s}{t^s}\big)^\frac{k+1}{2s}\notag\\
&\leq&\frac{\bar{B}_k}{t^\frac{k+1}{2}},
\end{eqnarray}
where we used \eqref{4.11} and \eqref{4.35}  and $Kt\leq a$ in the second inequality.

This completes the proof.
\end{proof}

Now we can prove Theorem \ref{cor1.1} with the aid of Theorem \ref{thm4.2}.
\\\\\textbf{Proof of Theorem \ref{cor1.1}:} We prove the result by contraction. If $T<\infty$ is a maximal existence time for a smooth solution $g(x,t)$ to the $(\alpha,\beta)$-Ricci-Yamabe flow \eqref{1.1} and
\[\sup_{M\times[0,T)}|Rm|<\infty.\]

From Theorem \ref{thm4.2}, we know that all covariant derivatives of the Riemannian curvature stay uniformly bounded along this flow. Hence, $g(\cdot,t)$ converges smoothly to a complete limit metric $g(\cdot,T)$. Then the flow \eqref{1.1} can be restarted from the new initial metric $g(\cdot,T)$ by the short-time existence result Theorem \ref{thm1.1}, which contradicts that the finite $T$ is a maximal existence time.

This proves this theorem.$\hfill\Box$
\section{Local curvature estimates and long time existence \Rmnum{2}} % (fold)
\label{sec:5}

Inspired by the works of Kotschwat-Munteanu-Wang \cite{KMW16} in the study of Ricci flow and Li-Yuan \cite{LY19} in the study of $\kappa$-LYZ flow, we prove that upper bounds for the norm of the Riemmanian curvature on an given geodesic ball of a smooth solution to the $(\alpha,\beta)$-Ricci-Yamabe flow \eqref{1.1} can be locally explicitly estimated in terms of its local $L^\infty$-norm on such a ball at the initial time and a uniform local bound for the Ricci curvature and the second-order derivative of the scalar curvature.

Let $g(x,t)$, $t\in[0,T']$ be a smooth solution to the $(\alpha,\beta)$-Ricci-Yamabe flow \eqref{1.1} with $0<\alpha\leq1$ and $\beta>-\frac{\alpha}{n-1}$ on an $n$-dimensional closed Riemannian manifold $M$.

Proceeding as in \cite{KMW16,LY19}, we prove a local curvature estimate without loss of generality by assuming that $M$ is a smooth complete manifold with
\begin{equation}\label{5.14}
|Ric|+|\nabla^2R|\leq L
\end{equation}
on $\Omega\times[0,T']$, where $\Omega=B_{g(0)}(x_0,\rho/\sqrt{L})\subset\subset M^n$ is a open geodesic ball under the initial metric for $x_0\in M$ and positive constants $L$ and $\rho$.

Choose the cutoff function
\begin{equation}\label{5.15}
\phi(x)=\big(\frac{\rho/\sqrt{L}-d_{g(0)}(x_0,x)}{\rho/\sqrt{L}}\big)_+
\end{equation}
which is  Lipschitz with support $\overline{B_{g_0}(x_0,\rho/\sqrt{L})}$.

Denote a universal positive constant that depends only on $n$ and $p$ by $c$ in this subsection.

In the following, we calculate some necessary differential inequalities.
\begin{proposition}\label{prp5.1}
For any $p\geq3$, we have
\begin{eqnarray}\label{5.12}
&&\frac{\alpha}{L}\int_M |\nabla Ric|^2|Rm|^{p-1}\phi^{2p}dVol_{g(t)}\notag\\
&\leq&-\frac{3}{2L}\cdot\frac{d}{dt}\int_M|Ric|^2|Rm|^{p-1}\phi^{2p} dVol_{g(t)}\notag\\
&&+\frac{(\alpha+|\beta|)^2cL}{\alpha}\int_M|\nabla Rm|^2|Rm|^{p-3}\phi^{2p} dVol_{g(t)}\notag\\
&&+\frac{(\alpha+|\beta|)^2cL}{\alpha}\int_M|Rm|^{p-1}|\nabla \phi|^2\phi^{2p-2} dVol_{g(t)}\notag\\
&&+(\alpha+|\beta|) cL\int_M |Rm|^p\phi^{2p}dVol_{g(t)}.
\end{eqnarray}
\end{proposition}

\begin{proof}
Direct computation on $\Omega\times[0,T]$ yields
\begin{eqnarray}\label{5.x}
\frac{1}{2}\frac{\partial}{\partial t}|Ric|^2&=&g^{ia}g^{kb}\big(\frac{\partial}{\partial t}R_{ik}\big)R_{ab}+g^{ie}g^{la}(2\alpha R_{le}+\beta Rg_{le})g^{kb}R_{ik}R_{ab}\notag\\
&=&\alpha g^{ia}g^{kb} R_{ab}\Delta R_{ik}+2\alpha g^{ia}g^{je}g^{kb}g^{ld}R_{ijkl}R_{ab}R_{de}\notag\\
&&-2\alpha g^{ja}g^{li}g^{bk}R_{jl}R_{ab}R_{ik}\notag\\
&&+\frac{\beta}{2}[(n-2)g^{ia}g^{kb}R_{ik}\nabla_a\nabla_bR+R\Delta R]\notag\\
&&+2\alpha g^{ja}g^{li}g^{bk}R_{jl}R_{ab}R_{ik}+\beta R|Ric|^2\notag\\
&\leq&\frac{\alpha}{2}\Delta|Ric|^2-\alpha|\nabla Ric|^2+\alpha cL^2|Rm|\notag\\
&&+|\beta|c\nabla^2R*Ric+\frac{\beta}{2}R\Delta R,
\end{eqnarray}
where we used \eqref{2.4} in the first equality and \eqref{2.6} in the second. Moreover, we used Cauchy's inequality and the fact of $|Ric|\leq L$ in the last inequality. By rearranging, we have
\begin{eqnarray}\label{5.2}
\alpha|\nabla Ric|^2&\leq&\frac{\alpha}{2}\Delta|Ric|^2-\frac{1}{2}\frac{\partial}{\partial t}|Ric|^2+\alpha cL^2|Rm|\notag\\
&&+|\beta|c\nabla^2R*Ric+\frac{\beta}{2}R\Delta R.
\end{eqnarray}

Similarly, we can derive from \eqref{2.1} that
\begin{equation}\label{5.3}
\alpha|\nabla Rm|^2\leq\frac{\alpha}{2}\Delta|Rm|^2-\frac{1}{2}\frac{\partial}{\partial t}|Rm|^2+(\alpha+|\beta|)c|Rm|^3+2|\beta|\cdot|\nabla^2R|\cdot|Ric|.
\end{equation}

Moreover,
\begin{eqnarray}\label{5.4}
\frac{\partial}{\partial t}|Rm|^2&=&\frac{\partial}{\partial t}(g^{ia}g^{jb}g^{ld}g_{kc}R_{ijl}^kR_{abd}^c)\notag\\
&\leq&(\alpha+|\beta|)c(\nabla^2Ric*Rm+\nabla^2R*Ric+Ric*Rm*Rm),
\end{eqnarray}
where we used \eqref{1.1}, \eqref{2.3} and \eqref{2.4}.

It follows from \eqref{5.2} that
\begin{eqnarray}\label{5.8}
&&\frac{\alpha}{L}\int_M |\nabla Ric|^2|Rm|^{p-1}\phi^{2p}dVol_{g(t)}\notag\\
&\leq&\frac{\alpha}{2L}\int_M (\Delta |Ric|^2)|Rm|^{p-1}\phi^{2p}dVol_{g(t)}\notag\\
&&+\frac{|\beta|}{L}\int_M(\nabla^2R*Ric)|Rm|^{p-1}\phi^{2p} dVol_{g(t)}\notag\\
&&+\frac{\beta}{2L}\int_MR\Delta R|Rm|^{p-1}\phi^{2p}dVol_{g(t)}\notag\\
&&-\frac{1}{2L}\int_M\big(\frac{\partial}{\partial t}|Ric|^2\big)|Rm|^{p-1}\phi^{2p} dVol_{g(t)}\notag\\
&&+\alpha cL\int_M |Rm|^p\phi^{2p}dVol_{g(t)}.
\end{eqnarray}

Integrating by parts and then using Cauchy's inequality yield
\begin{eqnarray}\label{5.9}
&&\frac{\alpha}{2L}\int_M (\Delta |Ric|^2)|Rm|^{p-1}\phi^{2p}dVol_{g(t)}\notag\\
&&+\frac{|\beta|}{L}\int_M(\nabla^2R*Ric)|Rm|^{p-1}\phi^{2p} dVol_{g(t)}\notag\\
&&+\frac{\beta}{2L}\int_MR\Delta R|Rm|^{p-1}\phi^{2p}dVol_{g(t)}\notag\\
&=&-\frac{\alpha}{2L}\int_M\langle\nabla|Ric|^2,\nabla|Rm|^{p-1}\rangle\phi^{2p} dVol_{g(t)}\notag\\
&&-\frac{\alpha}{2L}\int_M\langle\nabla|Ric|^2,\nabla(\phi^{2p})\rangle|Rm|^{p-1}\phi^{2p} dVol_{g(t)}\notag\\
&&+\frac{|\beta|}{L}\int_M(\nabla R*\nabla Ric)|Rm|^{p-1}\phi^{2p} dVol_{g(t)}\notag\\
&&+\frac{|\beta|}{L}\int_M(\nabla R*Ric*\nabla|Rm|^{p-1})\phi^{2p} dVol_{g(t)}\notag\\
&&+\frac{|\beta|}{L}\int_M(\nabla R*Ric*\nabla(\phi^{2p}))|Rm|^{p-1} dVol_{g(t)}\notag\\
&&-\frac{\beta}{2L}\int_M|\nabla R|^2|Rm|^{p-1}\phi^{2p}dVol_{g(t)}\notag\\
&&-\frac{\beta}{2L}\int_MR\langle\nabla R,\nabla|Rm|^{p-1}\rangle\phi^{2p}dVol_{g(t)}\notag\\
&&-\frac{\beta}{2L}\int_MR\langle\nabla R,\nabla(\phi^{2p})\rangle|Rm|^{p-1}dVol_{g(t)}\notag\\
&\leq&(\alpha+|\beta|)c\int_M |\nabla Ric|\cdot|\nabla Rm|\cdot|Rm|^{p-2}\phi^{2p}dVol_{g(t)}\notag\\
&&+(\alpha+|\beta|)c\int_M |\nabla Ric|\cdot|\nabla\phi|\cdot|Rm|^{p-1}\phi^{2p-1}dVol_{g(t)}\notag\\
&\leq&\frac{\alpha}{3L}\int_M |\nabla Ric|^2|Rm|^{p-1}\phi^{2p}dVol_{g(t)}\notag\\
&&+\frac{(\alpha+|\beta|)^2cL}{\alpha}\int_M |\nabla Rm|^2|Rm|^{p-3}\phi^{2p}dVol_{g(t)}\notag\\
&&+\frac{(\alpha+|\beta|)^2cL}{\alpha}\int_M |Rm|^{p-1}|\nabla\phi|^2\phi^{2p-2}dVol_{g(t)}.
\end{eqnarray}

%It follows from \eqref{5.4} and \eqref{5.5} that
Similar to \eqref{4.9}, we have
\begin{equation}\label{5.5}
\frac{\partial}{\partial t}dVol_{g(t)}=-(\alpha+\frac{n}{2}\beta)RdVol_{g(t)}\leq (\alpha+|\beta|)c|Ric|dVol_{g(t)}.
\end{equation}
Moreover,
\begin{equation}\label{5.6}
\frac{\partial}{\partial t}|Rm|^{p-1}=\frac{\partial}{\partial t}(|Rm|^2)^\frac{p-1}{2}=\frac{p-1}{2}|Rm|^{p-3}\frac{\partial}{\partial t}|Rm|^2
\end{equation}
for $p\geq3$.

Therefore,
\begin{eqnarray}\label{5.10}
&&-\frac{1}{2L}\int_M\big(\frac{\partial}{\partial t}|Ric|^2\big)|Rm|^{p-1}\phi^{2p} dVol_{g(t)}\notag\\
&=&-\frac{1}{2L}\cdot\frac{d}{dt}\int_M|Ric|^2|Rm|^{p-1}\phi^{2p} dVol_{g(t)}\notag\\
&&+\frac{p-1}{4L}\int_M |Ric|^2|Rm|^{p-3}\big(\frac{\partial}{\partial t}|Rm|^{2}\big)\phi^{2p}dVol_{g(t)}\notag\\
&&-\frac{1}{2L}\big(\alpha+\frac{n}{2}\beta\big)\int_MR|Ric|^2|Rm|^{p-1}\phi^{2p}dVol_{g(t)}\notag\\
&\leq&-\frac{1}{2L}\cdot\frac{d}{dt}\int_M|Ric|^2|Rm|^{p-1}\phi^{2p} dVol_{g(t)}\notag\\
&&+\frac{(\alpha+|\beta|)c}{L}\int_M |Ric|^2|Rm|^{p-3}(\nabla^2Ric*Rm+\nabla^2R*Ric)\phi^{2p}dVol_{g(t)}\notag\\
&&+(\alpha+|\beta|)cL\int_M|Rm|^{p}\phi^{2p}dVol_{g(t)},
\end{eqnarray}
where we used \eqref{5.4}. Integrating by parts and then using Cauchy's inequality, we obtain
\begin{eqnarray}\label{5.11}
&&\frac{(\alpha+|\beta|)c}{L}\int_M |Ric|^2|Rm|^{p-3}(\nabla^2Ric*Rm+\nabla^2R*Ric)\phi^{2p}dVol_{g(t)}\notag\\
&=&\frac{(\alpha+|\beta|)c}{L}\int_M \nabla|Ric|^2*(\nabla Ric*Rm+\nabla R*Ric)|Rm|^{p-3}\phi^{2p}dVol_{g(t)}\notag\\
&&+\frac{(\alpha+|\beta|)c}{L}\int_M |Ric|^2(\nabla Ric*Rm+\nabla R*Ric)*\nabla|Rm|^{p-3}\phi^{2p}dVol_{g(t)}\notag\\
&&+\frac{(\alpha+|\beta|)c}{L}\int_M |Ric|^2(\nabla Ric*\nabla Rm+\nabla R*\nabla Ric)|Rm|^{p-3}\phi^{2p}dVol_{g(t)}\notag\\
&&+\frac{(\alpha+|\beta|)c}{L}\int_M |Ric|^2(\nabla Ric*Rm+\nabla R*Ric)*\nabla(\phi^{2p})|Rm|^{p-3}dVol_{g(t)}\notag\\
&\leq&(\alpha+|\beta|)c\int_M |\nabla Ric|\cdot|\nabla Rm|\cdot|Rm|^{p-2}\phi^{2p}dVol_{g(t)}\notag\\
&&+(\alpha+|\beta|)c\int_M |\nabla Ric|\cdot|Rm|^{p-1}\cdot|\nabla\phi|\phi^{2p-1}dVol_{g(t)}\notag\\
&\leq&\frac{\alpha}{3L}\int_M|\nabla Ric|^2|Rm|^{p-1}\phi^{2p} dVol_{g(t)}\notag\\
&&+\frac{(\alpha+|\beta|)^2cL}{\alpha}\int_M|\nabla Rm|^2|Rm|^{p-3}\phi^{2p} dVol_{g(t)}\notag\\
&&+\frac{(\alpha+|\beta|)^2cL}{\alpha}\int_M|Rm|^{p-1}|\nabla \phi|^2\phi^{2p-2} dVol_{g(t)}.
\end{eqnarray}

Plugging \eqref{5.11} into \eqref{5.10} and then applying \eqref{5.9} and \eqref{5.10} to \eqref{5.8}, we conclude \eqref{5.12} immediately.
\end{proof}

\begin{lemma}\label{lem5.2}
Define
\begin{eqnarray*}
U(t)&:=&\int_M|Rm|^p\phi^{2p} dVol_{g(t)}+\frac{3}{2L}\int_M |Ric|^2|Rm|^{p-1}\phi^{2p}dVol_{g(t)}\notag\\
&&+\frac{(\alpha+|\beta|)^2cL}{\alpha^2}\int_M |Rm|^{p-1}\phi^{2p}dVol_{g(t)}.
\end{eqnarray*}
Then, for any $p\geq3$, we have
\begin{eqnarray}\label{5.1}
\frac{d}{dt}U(t)&\leq&(\alpha+|\beta|)\big(1+\frac{\beta^2}{\alpha^2}\big)cL\int_M|Rm|^p\phi^{2p} dVol_{g(t)}\notag\\
&&+\frac{(\alpha+|\beta|)^2}{\alpha}cL\int_M|Rm|^{p-1}|\nabla\phi|^2\phi^{2p-2} dVol_{g(t)}\notag\\
&&+\frac{(\alpha+|\beta|)^2|\beta|}{\alpha^2}cL^2\int_M|Rm|^{p-2}\phi^{2p}dVol_{g(t)}.
\end{eqnarray}
\end{lemma}

\begin{proof}
By direct computation, we have
\begin{eqnarray}\label{5.7}
&&\frac{d}{d t}\int_M |Rm|^p\phi^{2p}dVol_{g(t)}\notag\\
&=&\frac{p}{2}\int_M|Rm|^{p-2}\frac{\partial}{\partial t}|Rm|^2\phi^{2p} dVol_{g(t)}-(\alpha+\frac{n}{2}\beta)\int_M R|Rm|^p\phi^{2p}dVol_{g(t)}\notag\\
&\leq&(\alpha+|\beta|)c\big(\int_M |Rm|^{p-2}(\nabla^2Ric*Rm+\nabla^2R*Ric)\phi^{2p}dVol_{g(t)}\notag\\
&&+L\int_M |Rm|^p\phi^{2p}dVol_{g(t)}\big)\notag\\
&=&(\alpha+|\beta|)c\int_M \nabla|Rm|^{p-2}*(\nabla Ric*Rm+\nabla R*Ric)\phi^{2p}dVol_{g(t)}\notag\\
&&+(\alpha+|\beta|)c\int_M |Rm|^{p-2}(\nabla Ric*\nabla Rm+\nabla R*Ric)\phi^{2p}dVol_{g(t)}\notag\\
&&+(\alpha+|\beta|)c\int_M |Rm|^{p-2}(\nabla Ric* Rm+\nabla R*Ric)*\nabla(\phi^{2p})dVol_{g(t)}\notag\\
&&+(\alpha+|\beta|)cL\int_M |Rm|^{p}\phi^{2p}dVol_{g(t)}\notag\\
&\leq&(\alpha+|\beta|)c\int_M |\nabla Ric|\cdot|\nabla Rm|\cdot|Rm|^{p-2}\phi^{2p}dVol_{g(t)}\notag\\
&&+(\alpha+|\beta|)c\int_M |\nabla Ric|\cdot|Rm|^{p-1}\cdot|\nabla\phi|\phi^{2p-1}dVol_{g(t)}\notag\\
&&+(\alpha+|\beta|)cL\int_M |Rm|^{p}\phi^{2p}dVol_{g(t)}\notag\\
&\leq&\frac{\alpha}{L}\int_M |\nabla Ric|^2|Rm|^{p-1}\phi^{2p}dVol_{g(t)}\notag\\
&&+\frac{(\alpha+|\beta|)^2cL}{\alpha}\int_M|\nabla Rm|^2|Rm|^{p-3}\phi^{2p} dVol_{g(t)}\notag\\
&&+\frac{(\alpha+|\beta|)^2cL}{\alpha}\int_M|Rm|^{p-1}|\nabla\phi|^2\phi^{2p-2}dVol_{g(t)}\notag\\
&&+(\alpha+|\beta|)cL\int_M |Rm|^p\phi^{2p}dVol_{g(t)}\notag\\
&\leq&-\frac{3}{2L}\cdot\frac{d}{dt}\int_M|Ric|^2|Rm|^{p-1}\phi^{2p} dVol_{g(t)}\notag\\
&&+\frac{(\alpha+|\beta|)^2cL}{\alpha}\int_M|\nabla Rm|^2|Rm|^{p-3}\phi^{2p} dVol_{g(t)}\notag\\
&&+\frac{(\alpha+|\beta|)^2cL}{\alpha}\int_M|Rm|^{p-1}|\nabla \phi|^2\phi^{2p-2} dVol_{g(t)}\notag\\
&&+(\alpha+|\beta|) cL\int_M |Rm|^p\phi^{2p}dVol_{g(t)}.
\end{eqnarray}
where we used \eqref{5.5} and \eqref{5.6} in the first equality, \eqref{5.4} in the second equality,  Cauchy's inequality in the second and third inequalities and Proposition \ref{prp5.1} in the last.

Now we deal with the second term of RHS of \eqref{5.7}.
\begin{eqnarray}\label{5.13}
&&\int_M|\nabla Rm|^2|Rm|^{p-3}\phi^{2p} dVol_{g(t)}\notag\\
&\leq&\frac{1}{2}\int_M(\Delta|Rm|^2)|Rm|^{p-3}\phi^{2p} dVol_{g(t)}-\frac{1}{2\alpha}\int_M(\frac{\partial}{\partial t}|Rm|^2)|Rm|^{p-3}\phi^{2p} dVol_{g(t)}\notag\\
&&+\frac{(\alpha+|\beta|)c}{\alpha}\int_M|Rm|^{p}\phi^{2p} dVol_{g(t)}+\frac{4\beta}{\alpha}\int_M\langle\nabla^2R,Ric\rangle|Rm|^{p-3}\phi^{2p}dVol_{g(t)} \notag\\
&=&-\frac{1}{2}\int_M\langle\nabla|Rm|^2,\nabla|Rm|^{p-3}\rangle\phi^{2p} dVol_{g(t)}\notag\\
&&-\frac{1}{2}\int_M\langle\nabla|Rm|^2,\nabla \phi^{2p}\rangle|Rm|^{p-3} dVol_{g(t)}\notag\\
&&-\frac{1}{\alpha(p-1)}\int_M\frac{\partial}{\partial t}|Rm|^{p-1}\phi^{2p} dVol_{g(t)}+\frac{(\alpha+|\beta|)c}{\alpha}\int_M|Rm|^{p}\phi^{2p} dVol_{g(t)}\notag\\
&&+\frac{|\beta|}{\alpha}cL\int_M|Rm|^{p-2}\phi^{2p}dVol_{g(t)} \notag\\
&\leq&c\int_M |\nabla Rm|\cdot|\nabla\phi|\cdot|Rm|^{p-2}\cdot\phi^{2p-1}dVol_{g(t)}\notag\\
&&-\frac{1}{\alpha(p-1)}\cdot\frac{d}{d t}\int_M|Rm|^{p-1}\phi^{2p} dVol_{g(t)}\notag\\
&&+\frac{(\alpha+|\beta|)c}{\alpha}\int_M|Rm|^{p}\phi^{2p}dVol_{g(t)}+\frac{|\beta|}{\alpha}cL\int_M|Rm|^{p-2}\phi^{2p}dVol_{g(t)}\notag\\
&\leq&\frac{1}{2}\int_M|\nabla Rm|^2|Rm|^{p-3}\phi^{2p} dVol_{g(t)}+c\int_M|Rm|^{p-1}|\nabla \phi|^2\phi^{2p-2} dVol_{g(t)}\notag\\
&&-\frac{1}{\alpha(p-1)}\cdot\frac{d}{d t}\int_M|Rm|^{p-1}\phi^{2p} dVol_{g(t)}+\frac{(\alpha+|\beta|)c}{\alpha}\int_M|Rm|^{p}\phi^{2p}dVol_{g(t)}\notag\\
&&+\frac{|\beta|}{\alpha}cL\int_M|Rm|^{p-2}\phi^{2p}dVol_{g(t)},
\end{eqnarray}
where we used \eqref{5.3} in the first inequality, \eqref{5.5} in the second inequality and Cauchy's inequality in the last.

We conclude \eqref{5.1} by combining \eqref{5.7} and \eqref{5.13} and rearranging.
\end{proof}

In the rest of this section, we denote a universal constant depending only on $\alpha$, $\beta$, $L$ and $n$ by $C$. Then we have a $L^p$ estimate for the Riemannian curvature.
\begin{lemma}\label{lem5.3}
For any $x_0\in M^n$ and $\theta>1$, the Riemannian curvature of $g(t)$, $t\in[0,T']$, satisfies
\begin{eqnarray}\label{5.16}
&&\big(\int_{B_{g(0)}(x_0,\frac{\rho}{\theta\sqrt{L}})}|Rm|^p(x,t) dVol_{g(t)}\big)^\frac{1}{p}\notag\\
&\leq& \big(\frac{\theta}{\theta-1}\big)^{2p}Ce^{CLT'}\big[\big(\int_{B_{g(0)}(x_0,\frac{\rho}{\sqrt{L}})}|Rm|^p(x,0)dVol_{g(0)}\big)^\frac{1}{p}\notag\\
&&+\big(L^\frac{1}{2}+L+\frac{L}{\rho^{2}}\big)Vol^{\frac{1}{p}}_{g(0)}(B_{g(0)}(x_0,\frac{\rho}{\sqrt{L}}))\big].
\end{eqnarray}
\end{lemma}

\begin{proof}
It follows from \eqref{1.1} that on $B_{g(0)}(x_0,\rho/\sqrt{L})\times[0,T']$,
\[\frac{\partial}{\partial t}g\leq CLg,\]
i.e.,
\begin{equation}\label{5.17}
e^{-CLt}g(x,0)\leq g(x,t)\leq e^{CLt}g(x,0).
\end{equation}

Hence,
\begin{equation}\label{5.18}
Vol_{g(t)}(B_{g(0)}(x_0,\frac{\rho}{\sqrt{L}}))\leq e^{CLT'}Vol_{g(0)}(B_{g(0)}(x_0,\frac{\rho}{\sqrt{L}})),
\end{equation}
%and
%\begin{equation}\label{5.18x}
%Vol_{g(\tau)}(B_{g(0)}(x_0,\frac{\rho}{\sqrt{L}}))\leq e^{c(\alpha+|\beta|)LT}Vol_{g(0)}(B_{g(0)}(x_0,\frac{\rho}{\sqrt{L}}))
%\end{equation}
for all $0\leq t\leq T'$, and
\begin{equation}\label{5.19}
|\nabla\phi|_{g(t)}\leq e^{CLT'}|\nabla\phi|_{g(0)}\leq e^{CLT'}\cdot\frac{\sqrt{L}}{\rho}.
\end{equation}

Moreover,
\begin{eqnarray}\label{5.20}
&&\int_M |Rm|^{p-1}|\nabla\phi|^2\phi^{2p-2}dVol_{g(t)}\notag\\
&\leq&\frac{L}{\rho^2}\cdot e^{CLT'}\int_M |Rm|^{p-1}\phi^{2p-2}dVol_{g(t)}\notag\\
&\leq&\frac{L^p}{\rho^{2p}}\cdot e^{CLT'}Vol_{g(0)}(B_{g(0)}(x_0,\frac{\rho}{\sqrt{L}}))\notag\\
&&+\int_M |Rm|^{p}\phi^{2p}dVol_{g(t)},
\end{eqnarray}
where we used \eqref{5.19} in the first inequality, and Young's inequality and \eqref{5.18} in the last.

Using Young's inequality again, we can get
\begin{eqnarray}\label{5.20x}
&&L\int_M|Rm|^{p-2}\phi^{2p}dVol_{g(t)}\notag\\
&\leq&L^\frac{p}{2}\cdot e^{CLT'}Vol_{g(0)}(B_{g(0)}(x_0,\frac{\rho}{\sqrt{L}}))\notag\\
&&+\int_M |Rm|^{p}\phi^{2p}dVol_{g(t)},
\end{eqnarray}

Applying \eqref{5.20} and \eqref{5.20x} to \eqref{5.1}, we obtain
\begin{eqnarray*}
\frac{d}{dt}U(t)&\leq&CL\int_M|Rm|^p\phi^{2p}dVol_{g(t)}\notag\\
&&+C\big(L^{\frac{p}{2}+1}+\frac{L^{p+1}}{\rho^{2p}}\big)e^{CLT'}Vol_{g(0)}(B_{g(0)}(x_0,\frac{\rho}{\sqrt{L}}))\notag\\
&\leq&CLU(t)+C\big(L^{\frac{p}{2}+1}+\frac{L^{p+1}}{\rho^{2p}}\big)\notag\\
&&\times e^{CLT'}Vol_{g(0)}(B_{g(0)}(x_0,\frac{\rho}{\sqrt{L}})).
\end{eqnarray*}
This implies
\begin{eqnarray}\label{5.21}
&&\frac{d }{dt}\big(e^{-CLt}U(t)\big)\notag\\
&\leq&C\big(L^{\frac{p}{2}+1}+\frac{L^{p+1}}{\rho^{2p}}\big)e^{CL(T'-t)}Vol_{g(0)}(B_{g(0)}(x_0,\frac{\rho}{\sqrt{L}})).
\end{eqnarray}

For any $\tau\in(0,T']$, integrating \eqref{5.21} from $0$ to $\tau$ yields
\begin{equation}\label{5.22}
U(\tau)\leq e^{CLT'}\big(U(0)+C(L^\frac{p}{2}+\frac{L^p}{\rho^{2p}}) Vol_{g(0)}(B_{g(0)}(x_0,\frac{\rho}{\sqrt{L}}))\big).
\end{equation}

On the other hand, for any constant $\theta>1$, $\phi(x)\geq1-\frac{1}{\theta}$ for $d_{g(0)}(x_0,x)\leq\frac{\rho}{\theta\sqrt{L}}$. Therefore,
\begin{equation}\label{5.23}
\big(\frac{\theta-1}{\theta}\big)^{2p}\int_{B_{g(0)}(x_0,\frac{\rho}{\theta\sqrt{L}})}|Rm|^p(x,\tau) dVol_{g(\tau)}\leq\int_M|Rm|^p(x,\tau)\phi^{2p} dVol_{g(\tau)}\leq U(\tau).
\end{equation}

By definition, we have
\begin{eqnarray}\label{5.24}
U(0)&\leq&\int_M|Rm|^p(x,0)\phi^{2p} dVol_{g(0)}+\frac{3\sqrt{n}}{2}\int_M |Rm|^{p}(x,0)\phi^{2p}dVol_{g(0)}\notag\\
&&+CL\int_M |Rm|^{p-1}(x,0)\phi^{2p}dVol_{g(0)}\notag\\
&\leq&C\int_{B_{g(0)}(x_0,\frac{\rho}{\sqrt{L}})}|Rm|^p(x,0)\phi^{2p} dVol_{g(0)}\notag\\
&&+CL^pVol_{g(0)}(B_{g(0)}(x_0,\frac{\rho}{\sqrt{L}})),
\end{eqnarray}
where we used $|Ric|\leq L$ and $|Ric|\leq\sqrt{n}|Rm|$ in the first inequality and Young's inequality in the second.

It follows from \eqref{5.22}, \eqref{5.23} and \eqref{5.24} that
\begin{eqnarray}\label{5.25}
&&\int_{B_{g(0)}(x_0,\frac{\rho}{\theta\sqrt{L}})}|Rm|^p(x,\tau)\phi^{2p} dVol_{g(\tau)}\notag\\
&\leq&\big(\frac{\theta}{\theta-1}\big)^{2p}Ce^{CLT'}\big[\int_{B_{g(0)}(x_0,\frac{\rho}{\sqrt{L}})}|Rm|^p(x,0)dVol_{g(0)}\notag\\
&&+\big(L^\frac{p}{2}+L^p+\frac{L^{p}}{\rho^{2p}}\big)Vol_{g(0)}(B_{g(0)}(x_0,\frac{\rho}{\sqrt{L}}))\big]
\end{eqnarray}
for all $0<\tau\leq T'$ and $\theta>1$.

We conclude \eqref{5.16} from \eqref{5.25} by the arbitrary of $\tau$.
\end{proof}

Define
\[\dashint_B hdVol_{g(0)}=\frac{1}{Vol_{g(0)}(B)}\int_BhdVol_{g(0)}\]
to be the average of a smooth function $h\in M\times[0,T']$ on some geodesic ball $B$ with respect to $g(0)$. It is clear from \eqref{5.17} that
\begin{equation}\label{5.36}
e^{-CLT'}Vol_{g(0)}(B)\leq Vol_{g(t)}(B)\leq e^{CLT'}Vol_{g(0)}(B),
\end{equation}
for all $t\in(0,T']$.

Now we are ready to prove a local curvature estimate by applying De Giorgi-Nash-Moser iteration method (see e.g. \cite{PLi}). The proof follows from the work of Kotschwar, Munteanu and Wang \cite{KMW16} with minor modifications.

\begin{theorem}\label{thm5.1}
For any $x_0\in M$, there exist constants $a_1$ and $a_2$ depending only on $n$ so that
\begin{eqnarray}\label{5.26}
\sup_{B_{g(0)}(x_0,\frac{\rho}{4\sqrt{L}})\times[\frac{T'}{2},T']}|Rm|(x,t)&\leq& Ce^{C(LT'+\rho)}\big[1+\big(\frac{\Lambda_0}{L}\big)^{a_1}+\big(\frac{1}{LT'}+\frac{1}{\rho^2}\big)^{a_2}\big]\notag\\
&&\times\big(\Lambda_0+L^\frac{1}{2}+L+\frac{L}{\rho^2}\big),
\end{eqnarray}
where
\[\Lambda_0:=\sup_{B_{g(0)}(x_0,\frac{\rho}{\sqrt{L}})}|Rm|(x,0).\]
\end{theorem}

\begin{proof}
For $p\geq3$, it follows from \eqref{5.36} and Lemma \ref{lem5.3} that
\begin{eqnarray}\label{5.27}
&&\big(\dashint_{B_{g(0)}(x_0,\frac{\rho}{2\sqrt{L}})}|Rm|^p(x,t) dVol_{g(0)}\big)^\frac{1}{p}\notag\\
&\leq&e^{CLT'}\big(\frac{1}{Vol_{g(0)}(x_0,\frac{\rho}{2\sqrt{L}})}\int_{B_{g(0)}(x_0,\frac{\rho}{2\sqrt{L}})}|Rm|^p(x,t) dVol_{g(t)}\big)^\frac{1}{p}\notag\\
&\leq&Ce^{CLT'}\big(\Lambda_0+L^\frac{1}{2}+L+\frac{L}{\rho^2}\big)\big(\frac{Vol_{g(0)}(B_{g(0)}(x_0,\frac{\rho}{\sqrt{L}}))}{Vol_{g(0)}(B_{g(0)}(x_0,\frac{\rho}{2\sqrt{L}}))}\big)^\frac{1}{p}.
\end{eqnarray}

Since the Ricci curvature of $g(0)$ is lower bounded, by applying the Bishop-Gromov volume comparison theorem (see e.g. \cite{PLi}) to \eqref{5.27}, we have
\begin{eqnarray}\label{5.28}
&&\big(\dashint_{B_{g(0)}(x_0,\frac{\rho}{2\sqrt{L}})}|Rm|^p(x,t) dVol_{g(0)}\big)^\frac{1}{p}\notag\\
&\leq& Ce^{C(T'+\frac{\rho}{\sqrt{L}})}\big(\Lambda_0+L^\frac{1}{2}+L+\frac{L}{\rho^2}\big)
\end{eqnarray}
for $p\geq3$.

%Note that
%\[\Delta|Rm|^2=2|\nabla|Rm||^2+2|Rm|\Delta|Rm|\]
Note that $0<\alpha\leq1$, using Kato's inequality and \eqref{5.3}, we obtain
\begin{eqnarray*}
\alpha|\nabla |Rm||^2&\leq&\alpha|\nabla Rm|^2\notag\\
&\leq&\frac{\alpha}{2}\Delta|Rm|^2-\frac{1}{2}\frac{\partial}{\partial t}|Rm|^2+C|Rm|^3+C|\nabla^2R|\cdot|Rm|^2\notag\\
&\leq&|\nabla|Rm||^2+\alpha|Rm|\Delta|Rm|-|Rm|\frac{\partial}{\partial t}|Rm|+C|Rm|^3+CL|Rm|^2,
\end{eqnarray*}
i.e.,
\[\big(\frac{\partial}{\partial t}-\alpha\Delta\big)|Rm|\leq C(|Rm|+L)|Rm|\]
weakly on $M\times[0,T']$.

As in \cite{KMW16,LZ21ax}, we write $f=|Rm|+L$ and $u=|Rm|$, then
\begin{equation}\label{5.29}
\big(\frac{\partial}{\partial t}-\alpha\Delta\big)u\leq Cfu
\end{equation}
weakly on $M\times[0,T']$. Fixing $a\geq1$ and integrating \eqref{5.29} on $B_{g(0)}(x_0,\frac{\rho}{2\sqrt{L}})$ give that
\begin{eqnarray}\label{5.30}
&&\int_{B_{g(0)}(x_0,\frac{\rho}{2\sqrt{L}})}\varphi^2u^{2a-1}\frac{\partial u}{\partial t}dVol_{g(t)}-\alpha\int_{B_{g(0)}(x_0,\frac{\rho}{2\sqrt{L}})}\varphi^2u^{2a-1}\Delta u dVol_{g(t)}\notag\\
&\leq& C\int_{B_{g(0)}(x_0,\frac{\rho}{2\sqrt{L}})}\varphi^2u^{2a}fdVol_{g(t)}
\end{eqnarray}
for any $t\in[0,T']$ and nonnegative Lipschitz function $\varphi(x)$ with compact support in $B_{g(0)}(x_0,\frac{\rho}{2\sqrt{L}})$.

By \eqref{4.9}, we know that
\begin{eqnarray}\label{5.31}
\int_{B_{g(0)}(x_0,\frac{\rho}{2\sqrt{L}})}\varphi^2u^{2a-1}\frac{\partial u}{\partial t}dVol_{g(t)}&\geq&\frac{1}{2a}\cdot\frac{d}{dt}\int_{B_{g(0)}(x_0,\frac{\rho}{2\sqrt{L}})}\varphi^2u^{2a}dVol_{g(t)}\notag\\
&&-C\int_{B_{g(0)}(x_0,\frac{\rho}{2\sqrt{L}})}\varphi^2u^{2a}fdVol_{g(t)}.
\end{eqnarray}

Integrating by parts yields
\begin{eqnarray}\label{5.32}
&&-\int_{B_{g(0)}(x_0,\frac{\rho}{2\sqrt{L}})}\varphi^2u^{2a-1}\Delta u dVol_{g(t)}\notag\\
&=& 2\int_{B_{g(0)}(x_0,\frac{\rho}{2\sqrt{L}})}\varphi^2u^{2a-1}\langle \nabla u,\nabla\varphi\rangle dVol_{g(t)}\notag\\
&&+(2a-1)\int_{B_{g(0)}(x_0,\frac{\rho}{2\sqrt{L}})}\varphi^2u^{2a-2}|\nabla u|^2 dVol_{g(t)}\notag\\
&\geq&2\int_{B_{g(0)}(x_0,\frac{\rho}{2\sqrt{L}})}\varphi^2u^{2a-1}\langle \nabla u,\nabla\varphi\rangle dVol_{g(t)}\notag\\
&&+a\int_{B_{g(0)}(x_0,\frac{\rho}{2\sqrt{L}})}\varphi^2u^{2a-2}|\nabla u|^2 dVol_{g(t)}\notag\\
&=&\frac{1}{a}\int_{B_{g(0)}(x_0,\frac{\rho}{2\sqrt{L}})}|\nabla(\varphi u^a)|^2dVol_{g(t)}\notag\\
&&-\frac{1}{a}\int_{B_{g(0)}(x_0,\frac{\rho}{2\sqrt{L}})}|\nabla\varphi|^2u^{2a}dVol_{g(t)}.
\end{eqnarray}

Plugging \eqref{5.31} and \eqref{5.32} into \eqref{5.30}, we have
\begin{eqnarray}\label{5.33}
&&\frac{1}{2}\cdot\frac{d}{dt}\int_{B_{g(0)}(x_0,\frac{\rho}{2\sqrt{L}})}\varphi^2u^{2a}dVol_{g(t)}+\alpha\int_{B_{g(0)}(x_0,\frac{\rho}{2\sqrt{L}})}|\nabla(\varphi u^a)|^2dVol_{g(t)}\notag\\
&\leq&Ca\int_{B_{g(0)}(x_0,\frac{\rho}{2\sqrt{L}})}\varphi^2u^{2a}fdVol_{g(t)}+\int_{B_{g(0)}(x_0,\frac{\rho}{2\sqrt{L}})}|\nabla\varphi|^2u^{2a}dVol_{g(t)}.
\end{eqnarray}

For $0<s<s+v<T'$, multiply \eqref{5.33} with the Lipschitz function $\psi(t)$ that is as same as (3.6) in \cite{KMW16},
\[\psi(t):=
\begin{cases}
0 & \text{for}\ 0\leq t\leq s\\
\frac{t-s}{v} &\text{for}\ s<t\leq s+v\\
1 &\text{for}\ s+v<t\leq T
\end{cases}
\]
and obtain
\begin{eqnarray}\label{5.34}
&&\frac{1}{2}\cdot\frac{d}{dt}\big(\psi^2\int_{B_{g(0)}(x_0,\frac{\rho}{2\sqrt{L}})}\varphi^2u^{2a}dVol_{g(t)}\big)\notag\\
&&+\alpha\psi^2\int_{B_{g(0)}(x_0,\frac{\rho}{2\sqrt{L}})}|\nabla(\varphi u^a)|^2dVol_{g(t)}\notag\\
&\leq&Ca\psi^2\int_{B_{g(0)}(x_0,\frac{\rho}{2\sqrt{L}})}\varphi^2u^{2a}fdVol_{g(t)}\notag\\
&&+\psi^2\int_{B_{g(0)}(x_0,\frac{\rho}{2\sqrt{L}})}|\nabla\varphi|^2u^{2a}dVol_{g(t)}\notag\\
&&+\psi\psi'\int_{B_{g(0)}(x_0,\frac{\rho}{2\sqrt{L}})}\varphi^2u^{2a}dVol_{g(t)}.
\end{eqnarray}

For $\tau\in(0,T']$, integrating \eqref{5.34} from $0$ to $\tau$ and then using \eqref{5.36} yield
\begin{eqnarray}\label{5.35}
&&\frac{1}{2}\psi^2(\tau)\dashint_{B_{g(0)}(x_0,\frac{\rho}{2\sqrt{L}})}\varphi^2u^{2a}(\tau)dVol_{g(0)}\notag\\
&&+\alpha\int_0^\tau\psi^2\dashint_{B_{g(0)}(x_0,\frac{\rho}{2\sqrt{L}})}|\nabla(\varphi u^a)|^2dVol_{g(0)}\notag\\
&\leq&Ca e^{CLT}\int_0^\tau\psi^2\dashint_{B_{g(0)}(x_0,\frac{\rho}{2\sqrt{L}})}\varphi^2u^{2a}fdVol_{g(0)}\notag\\
&&+e^{CLT}\int_0^\tau\psi^2\dashint_{B_{g(0)}(x_0,\frac{\rho}{2\sqrt{L}})}|\nabla\varphi|^2u^{2a}dVol_{g(0)}\notag\\
&&+e^{CLT}\int_0^\tau\psi\psi'\dashint_{B_{g(0)}(x_0,\frac{\rho}{2\sqrt{L}})}\varphi^2u^{2a}dVol_{g(0)}.
\end{eqnarray}

The condition of $|Ric|\leq L$ on $B_{g(0)}(x_0,\frac{\rho}{2\sqrt{L}})$ together with \eqref{5.17} give that there exists a Sobolev inequality of the form
\begin{equation}\label{5.37}
\frac{1}{C}\cdot\frac{e^{-C\rho}}{\rho^2}\big(\dashint_{B_{g(0)}(x_0,\frac{\rho}{2\sqrt{L}})}(\varphi u^a(t))^{2\mu}dVol_{g(0)}\big)^\frac{1}{\mu}\leq\dashint_{B_{g(0)}(x_0,\frac{\rho}{2\sqrt{L}})}|\nabla(\varphi u^a(t))|_{g(0)}^2dVol_{g(0)}
\end{equation}
for some constant $\mu=\mu(n)\leq\frac{n}{n-2}$ depending only on $n$.

Applying \eqref{5.37} to \eqref{5.35}, we have
\begin{eqnarray}\label{5.38}
&&\frac{1}{2}\psi^2(\tau)\dashint_{B_{g(0)}(x_0,\frac{\rho}{2\sqrt{L}})}\varphi^2u^{2a}(\tau)dVol_{g(0)}\notag\\
&&+\frac{1}{C}\frac{e^{-C\rho}}{\rho^2}\int_0^\tau\big(\dashint_{B_{g(0)}(x_0,\frac{\rho}{2\sqrt{L}})}(\varphi u^a(t))^{2\mu}dVol_{g(0)}\big)^\frac{1}{\mu}\notag\\
&\leq&Cae^{CLT}\int_0^\tau\psi^2\dashint_{B_{g(0)}(x_0,\frac{\rho}{2\sqrt{L}})}\varphi^2u^{2a}fdVol_{g(0)}\notag\\
&&+e^{CLT}\int_0^\tau\psi^2\dashint_{B_{g(0)}(x_0,\frac{\rho}{2\sqrt{L}})}|\nabla\varphi|^2u^{2a}dVol_{g(0)}\notag\\
&&+e^{CLT}\int_0^\tau\psi\psi'\dashint_{B_{g(0)}(x_0,\frac{\rho}{2\sqrt{L}})}\varphi^2u^{2a}dVol_{g(0)},
\end{eqnarray}
which is exactly as same as (3.9) in \cite{KMW16}. Following the iteration arguments of (3.10) to (3.11) in \cite{KMW16} for the rest steps, together with \eqref{5.28}, we conclude \eqref{5.26}.
\end{proof}

Finally, we give a proof of Theorem \ref{cor1.2}.
\\\\\textbf{Proof of Theorem \ref{cor1.2}:}
Since the Ricci curvature is uniformly bounded, it follows from Theorem \ref{thm5.1} and the compactness of $M$ that the Riemannian curvature of $g(t)$ is uniformly bounded along the $(\alpha,\beta)$-Ricci-Yamabe flow \ref{1.1}. Therefore, this flow can be extended smoothly by Theorem \ref{cor1.2}.$\hfill\Box$

\section{Curvature pinching estimates and long time existence \Rmnum{3}}
\label{sec:6}
In this section, we derive curvature pinching estimates on the traceless Ricci curvature $\overset\circ{Ric}:=Ric-\frac{1}{n}Rg$ under an $n$-dimensional $(n\geq3)$ $(\alpha,\beta)$-Ricci-Yamabe flow \eqref{1.1}, which generalize Cao's result \cite{cao11}. As an application, we show finite-time singularity behaviors of \eqref{1.1} in terms of scalar curvature and Weyl tensor.

Choose $b=2\max_M|R(\cdot,0)|+1$, it follows from Proposition \ref{prp2.1} that
\begin{equation}\label{6.12}
b\geq1,\ R+b\geq1\ \text{and}\ (R+b)^2\geq R^2.
\end{equation}
along the $(\alpha,\beta)$-Ricci-Yamabe flow \eqref{1.1} with $\alpha>0$ and $\beta>-\frac{\alpha}{n-1}$.

As in \cite{cao11,lirhf18}, we define
\begin{eqnarray}\label{6.1}
f&:=&\frac{|\overset\circ{Ric}|^2}{(R+b)^2}\notag\\
&=&\frac{|Ric+\frac{b}{n}g-\frac{R+b}{n}g|^2}{(R+b)^2}\notag\\
&=&\frac{|Ric|^2}{(R+b)^2}+\frac{2b R}{n(R+b)^2}+\frac{b^2}{n(R+b)^2}-\frac{1}{n}\notag\\
&=&\frac{|Ric|^2}{(R+b)^2}+\frac{R^2}{n(R+b)^2}.
\end{eqnarray}

For any $(x,t)\in M^n\times[0,T)$, we calculate in a normal coordinate system centered at $x$.
\begin{proposition}\label{prp6.1}
Under $(\alpha,\beta)$-Ricci-Yamabe flow \eqref{1.1}, we have
\begin{eqnarray}\label{6.2}
&&\frac{1}{(R+b)^2}\frac{\partial}{\partial t}|Ric+\frac{b}{n}g|^2\notag\\
&=&\alpha\Delta f+2\alpha\frac{\langle\nabla f,\nabla R\rangle}{R+b}+\frac{2\alpha b^2|\nabla R|^2}{n(R+b)^4}\notag\\
&&+[\frac{2\alpha\Delta R}{R+b}+\frac{4\alpha(2n-3)}{n(n-1)}R+\beta R+\frac{4b\alpha}{n}]f\notag\\
&&+[\frac{2\alpha R}{n}+\beta R+\frac{2b}{n}(\beta(n-1)-\alpha(n-2))]\frac{\Delta R}{(R+b)^2}\notag\\
&&+[\frac{4\alpha(n-8)}{n(n-2)}R+\beta R+2b\beta+\frac{4b\alpha}{n}]\frac{R^2}{n(R+b)^2}\notag\\
&&+4\alpha \frac{W(\overset\circ{Ric},\overset\circ{Ric})}{(R+b)^2}-\frac{8\alpha}{n-2}\frac{tr\overset\circ{Ric^3}}{(R+b)^2}+(n-2)\beta\frac{\langle Ric,\nabla^2R\rangle}{(R+b)^2}\notag\\
&&-2\alpha \frac{|(R+b)\nabla_iR_{jk}-\nabla_iRR_{jk}|^2}{(R+b)^4},
\end{eqnarray}
where $W(\overset\circ{Ric},\overset\circ{Ric}):=W_{ijkl}\overset\circ{R}_{ik}\overset\circ{R}_{jl}$ and $tr\overset\circ{Ric^3}=R_{ij}R_{jk}R_{ki}$.
\end{proposition}

\begin{proof}
Similar as the computation in \eqref{5.x}, it follows from \eqref{2.6} and \eqref{2.7} that
\begin{eqnarray}\label{6.3}
&&\frac{\partial}{\partial t}|Ric+\frac{b}{n}g|^2\notag\\
&=&\frac{\partial}{\partial t}|Ric|^2+\frac{2b}{n}\frac{\partial R}{\partial t}\notag\\
&=&\alpha\Delta|Ric|^2-2\alpha|\nabla Ric|^2+4\alpha Rm(Ric,Ric)\notag\\
&&+\beta R|Ric|^2+(n-2)\beta\langle Ric,\nabla^2R\rangle+\beta R\Delta R\notag\\
&&+\frac{2b}{n}[\beta(n-1)+\alpha]\Delta R+\frac{4b\alpha}{n}|Ric|^2+\frac{2b\beta}{n}R^2.
\end{eqnarray}

By \eqref{6.1}, we have
\begin{eqnarray}\label{6.4}
\langle\nabla f,\nabla R\rangle&=&\frac{\langle\nabla|Ric|^2,\nabla R\rangle}{(R+b)^2}+\frac{2b|\nabla R|^2}{n(R+b)^2}\notag\\
&&-\frac{2f|\nabla R|^2}{R+b}-\frac{2|\nabla R|^2}{n(R+b)}
\end{eqnarray}
and
\begin{eqnarray}\label{6.5}
\Delta f&=&div\big(\frac{\nabla|Ric|^2}{(R+b)^2}+\frac{2b\nabla R}{n(R+b)^2}-\frac{2f\nabla R}{R+b}-\frac{2\nabla R}{n(R+b)}\big)\notag\\
&=&\frac{\Delta|Ric|^2}{(R+b)^2}-\frac{2\langle\nabla|Ric|^2,\nabla R\rangle}{(R+b)^3}+\frac{2b\Delta R}{n(R+b)^2}\notag\\
&&-\frac{4b|\nabla R|^2}{n(R+b)^3}-\frac{2\langle \nabla f,\nabla R\rangle}{R+b}-\frac{2f\Delta R}{R+b}\notag\\
&&+\frac{2f|\nabla R|^2}{(R+b)^2}-\frac{2\Delta R}{n(R+b)}+\frac{2|\nabla R|^2}{n(R+b)^2}\notag\\
&=&\frac{\Delta|Ric|^2}{(R+b)^2}-\frac{4\langle \nabla f,\nabla R\rangle}{R+b}-\frac{2f|\nabla R|^2}{(R+b)^2}-\frac{2|\nabla R|^2}{n(R+b)^2}\notag\\
&&+\frac{2b\Delta R}{(R+b)^2}-\frac{2f\Delta R}{R+b}-\frac{2\Delta R}{n(R+b)}.
\end{eqnarray}

Note that
\begin{eqnarray}\label{6.6}
|\nabla Ric|^2&=&\frac{|(R+b)\nabla_iR_{jk}-\nabla_iRR_{jk}|^2}{(R+b)^2}+\frac{\langle\nabla|Ric|^2,\nabla R\rangle}{R+b}\notag\\
&&-\frac{|\nabla R|^2|Ric|^2}{(R+b)^2}\notag\\
&=&\frac{|(R+b)\nabla_iR_{jk}-\nabla_iRR_{jk}|^2}{(R+b)^2}+(R+b)\langle\nabla f,\nabla R\rangle\notag\\
&&-\frac{2b|\nabla R|^2}{n(R+b)}+f|\nabla R|^2+\frac{1}{n}|\nabla R|^2\notag\\
&&+\frac{2bR|\nabla R|^2}{n(R+b)^2}+\frac{b^2|\nabla R|^2}{n(R+b)^2},
\end{eqnarray}
where we used \eqref{6.4} and \eqref{6.1}.

By the definition of the Weyl tensor, we have
\begin{eqnarray*}
Rm(Ric,Ric)&=&R_{ijkl}R_{ik}R_{jl}\notag\\
&=&W_{ijkl}R_{ik}R_{jl}+\frac{2}{n-2}R|Ric|^2\notag\\
&&-\frac{1}{n-2}(R_{jk}R_{ik}R_{ij}+R_{il}R_{ik}R_{kl})\notag\\
&&-\frac{1}{(n-1)(n-2)}R(R^2-|Ric|^2)\notag\\
&=&W(Ric,Ric)+\frac{1}{n-2}\big(\frac{2n-1}{n-1}|Ric|^2R-2trRic^3-\frac{R^3}{n-1}\big).
\end{eqnarray*}

Since the Weyl tensor is trace-less and
\begin{eqnarray*}
trRic^3&=&R_{ij}R_{jk}R_{ki}\notag\\
&=&(\overset\circ{R}_{ij}+\frac{R}{n}g_{ij})(\overset\circ{R}_{jk}+\frac{R}{n}g_{jk})(\overset\circ{R}_{ki}+\frac{R}{n}g_{ki})\notag\\
&=&(\overset\circ{R}_{ij}\overset\circ{R}_{jk}+\frac{2R}{n}\overset\circ{R}_{ik}+\frac{R^2}{n^2}g_{ik})(\overset\circ{R}_{ki}+\frac{R}{n}g_{ki})\notag\\
&=&tr\overset\circ{Ric^3}+\frac{3}{n}R|\overset\circ{Ric}|^2+\frac{R^3}{n^2},
\end{eqnarray*}
we get
\begin{eqnarray}\label{6.7}
Rm(Ric,Ric)&=&W(\overset\circ{Ric},\overset\circ{Ric})+\frac{2n-3}{n(n-1)}|Ric|^2R\notag\\
&&-\frac{2}{n-2}tr\overset\circ{Ric^3}-\frac{n^2+2n-2}{n^2(n-1)(n-2)}R^3\notag\\
&=&W(\overset\circ{Ric},\overset\circ{Ric})+\frac{2n-3}{n(n-1)}(R+b)^2fR\notag\\
&&-\frac{2}{n-2}tr\overset\circ{Ric^3}+\frac{n-8}{n^2(n-2)}R^3.
\end{eqnarray}

Applying \eqref{6.1}, \eqref{6.4}, \eqref{6.5}, \eqref{6.6}, \eqref{6.7} to \eqref{6.3}, we obtain
\begin{eqnarray}\label{6.8}
&&\frac{\partial}{\partial t}|Ric+\frac{b}{n}g|^2\notag\\
&=&\alpha(R+b)^2\Delta f+4\alpha(R+b)\langle\nabla f,\nabla R\rangle+2\alpha f|\nabla R|^2+\frac{2\alpha}{n}|\nabla R|^2\notag\\
&&-2\alpha b\Delta R+2\alpha f(R+b)\Delta R+\frac{2\alpha}{n}(R+b)\Delta R\notag\\
&&-2\alpha \frac{|(R+b)\nabla_iR_{jk}-\nabla_iRR_{jk}|^2}{(R+b)^2}-2\alpha(R+b)\langle\nabla f,\nabla R\rangle\notag\\
&&+\frac{4\alpha b|\nabla R|^2}{n(R+b)}-2\alpha f|\nabla R|^2-\frac{2\alpha}{n}|\nabla R|^2-\frac{4\alpha bR|\nabla R|^2}{n(R+b)^2}\notag\\
&&-\frac{2\alpha b^2|\nabla R|^2}{n(R+b)^2}+4\alpha W(\overset\circ{Ric},\overset\circ{Ric})+\frac{4\alpha(2n-3)}{n(n-1)}(R+b)^2fR\notag\\
&&-\frac{8\alpha}{n-2}tr\overset\circ{Ric^3}+\frac{4\alpha(n-8)}{n^2(n-2)}R^3+\beta(R+b)^2Rf+\frac{\beta}{n}R^3\notag\\
&&+\frac{2b\beta}{n}R^2+(n-2)\beta\langle Ric,\nabla^2R\rangle+\beta R\Delta R\notag\\
&&+\frac{2b}{n}[\beta(n-1)+\alpha]\Delta R+\frac{4b\alpha}{n}[(R+b)^2f+\frac{R^2}{n}]\notag\\
&=&-2\alpha \frac{|(R+b)\nabla_iR_{jk}-\nabla_iRR_{jk}|^2}{(R+b)^2}\notag\\
&&+\alpha(R+b)^2\Delta f+2\alpha(R+b)\langle\nabla f,\nabla R\rangle+\frac{2\alpha b^2|\nabla R|^2}{n(R+b)^2}\notag\\
&&+(R+b)^2[\frac{2\alpha\Delta R}{R+b}+\frac{4\alpha(2n-3)}{n(n-1)}R+\beta R+\frac{4b\alpha}{n}]f\notag\\
&&+[\frac{2\alpha R}{n}+\beta R+\frac{2b}{n}(\beta(n-1)-\alpha(n-2))]\Delta R\notag\\
&&+[\frac{4\alpha(n-8)}{n(n-2)}R+\beta R+2b\beta+\frac{4b\alpha}{n}]\frac{R^2}{n}\notag\\
&&+4\alpha W(\overset\circ{Ric},\overset\circ{Ric})-\frac{8\alpha}{n-2}tr\overset\circ{Ric^3}+(n-2)\beta\langle Ric,\nabla^2R\rangle.
\end{eqnarray}
Hence, \eqref{6.2} follows immediately.
\end{proof}

\begin{proposition}\label{prp6.2}
Under $(\alpha,\beta)$-Ricci-Yamabe flow \eqref{1.1}, we have
\begin{eqnarray}\label{6.9}
&&|Ric+\frac{b}{n}g|^2\frac{\partial}{\partial t}\frac{1}{(R+b)^2}\notag\\
&=&-4\alpha(R+b)f^2-\big[\frac{2[\beta(n-1)+\alpha]\Delta R}{R+b}\notag\\
&&+\frac{4}{n}\alpha(R+b)+2\big(\frac{2\alpha}{n}+\beta\big)\frac{R^2}{R+b}\big]f\notag\\
&&-\frac{2[\beta(n-1)+\alpha]\Delta R}{n(R+b)}-2\big(\frac{2\alpha}{n}+\beta\big)\frac{R^2}{n(R+b)}.
\end{eqnarray}
\end{proposition}

\begin{proof}
By direct computation, we obtain
\begin{eqnarray}\label{6.10}
\frac{\partial}{\partial t}\frac{1}{(R+b)^2}&=&-\frac{2}{(R+b)^3}\frac{\partial R}{\partial t}\notag\\
&=&-\frac{2[\beta(n-1)+\alpha]\Delta R}{(R+b)^3}-\frac{4\alpha|Ric|^2}{(R+b)^3}-\frac{2\beta R^2}{(R+b)^3}\notag\\
&=&-\frac{2[\beta(n-1)+\alpha]\Delta R}{(R+b)^3}-\frac{4\alpha}{R+b}f\notag\\
&&-2\big(\frac{2\alpha}{n}+\beta\big)\frac{R^2}{(R+b)^3},
\end{eqnarray}
where we used \eqref{2.7} in the second equality and \eqref{6.1} in the third.

It follows from \eqref{6.1} and \eqref{6.10} that
\begin{eqnarray*}
&&|Ric+\frac{b}{n}g|^2\frac{\partial}{\partial t}\frac{1}{(R+b)^2}\notag\\
&=&(f+\frac{1}{n})\big[-\frac{2[\beta(n-1)+\alpha]\Delta R}{R+b}\notag\\
&&-4\alpha(R+b)f-2\big(\frac{2\alpha}{n}+\beta\big)\frac{R^2}{R+b}\big]\notag\\
&=&-4\alpha(R+b)f^2-\big[\frac{2[\beta(n-1)+\alpha]\Delta R}{R+b}\notag\\
&&+\frac{4}{n}\alpha(R+b)+2\big(\frac{2\alpha}{n}+\beta\big)\frac{R^2}{R+b}\big]f\notag\\
&&-\frac{2[\beta(n-1)+\alpha]\Delta R}{n(R+b)}-2\big(\frac{2\alpha}{n}+\beta\big)\frac{R^2}{n(R+b)}.
\end{eqnarray*}
\end{proof}

Then we derive a inequality that are essential in the curvature pinching estimate.
\begin{lemma}\label{lem6.1}
Under $(\alpha,\beta)$-Ricci-Yamabe flow with $\alpha>0$, then there exists a uniform positive constant $\bar{c}$ depending only on $n$, $\alpha$ and $\beta$ so that
\begin{eqnarray}\label{6.11}
\frac{\partial}{\partial t}f&\leq&\alpha\Delta f+2\alpha\frac{\langle\nabla f,\nabla R\rangle}{R+b}-\alpha f^2\notag\\
&&\bar{c}\big(\frac{b^2|\nabla R|^2+|\nabla^2R|^2+|W|^2}{(R+b)^2}\big)+\bar{c}(R+b)^2.
\end{eqnarray}
\end{lemma}

\begin{proof}
By Proposition \ref{prp6.2}, Proposition \ref{prp6.2}, \eqref{6.1} and Cauchy's inequality, we have
\begin{eqnarray}\label{6.14}
\frac{\partial}{\partial t}f&=&\frac{1}{(R+b)^2}\frac{\partial}{\partial t}|Ric+\frac{b}{n}g|^2+|Ric+\frac{b}{n}g|^2\frac{\partial}{\partial t}\frac{1}{(R+b)^2}\notag\\
&\leq&\alpha\Delta f+2\alpha\frac{\langle\nabla f,\nabla R\rangle}{R+b}+\frac{2\alpha b^2|\nabla R|^2}{n(R+b)^4}\notag\\
&&+\bar{c}\big(\frac{|\nabla^2R|}{R+b}+R+b\big)f+\frac{\bar{c}(|\nabla^2R|+R^2)}{R+b}\notag\\
&&+\bar{c}|W|f+\bar{c}(R+b)f^\frac{3}{2}+\bar{c}(f+\frac{R^2+|\nabla^2R|^2}{(R+b)^2})\notag\\
&&-4\alpha(R+b)f^2+\bar{c}\big(\frac{|\nabla^2R|+R^2}{R+b}+R+b\big)f\notag\\
\end{eqnarray}

Applying \eqref{6.12} and Young's inequality to \eqref{6.14}, we obtain
\begin{eqnarray*}
\frac{\partial}{\partial t}f&\leq&\alpha\Delta f+2\alpha\frac{\langle\nabla f,\nabla R\rangle}{R+b}-4\alpha(R+b)f^2\notag\\
&&+\bar{c}\big(|\nabla^2R|+R+b+|W|\big)f+\bar{c}(R+b)f^\frac{3}{2}\notag\\
&&+\frac{\bar{c}b^2|\nabla R|^2}{(R+b)^2}+\frac{\bar{c}|\nabla^2 R|^2}{(R+b)^2}+\bar{c}(R+b)\notag\\
&\leq&\alpha\Delta f+2\alpha\frac{\langle\nabla f,\nabla R\rangle}{R+b}-\alpha f^2\notag\\
&&\bar{c}\big(\frac{b|\nabla R|+|\nabla^2R|+|W|}{R+b}\big)^2+\bar{c}(R+b)^2.
\end{eqnarray*}
\end{proof}

Applying maximum principle to \eqref{6.11}, we obtain the curvature pinching estimate immediately.
\begin{theorem}\label{thm6.1}
Let $g(x,t)$, $0\leq t<T<+\infty$, be a smooth solution to the $(\alpha,\beta)$-Ricci-Yamabe flow \eqref{1.1} with $\alpha>0$ and $\beta>-\frac{\alpha}{n-1}$on an $n$-dimensional $(n\geq3)$ closed  Riemannian manifold $M$. Define $b:=2\max_M|R(\cdot,0)|+1$. Then there exists a uniform positive constant $\tilde{C}$ depending only on $g(\cdot,0)$, $n$, $\alpha$ and $\beta$ so that
\begin{equation}\label{6.13}
\frac{|\overset\circ{Ric}|^2}{(R+b)^2}(x,t)\leq\bar{C}+\bar{C}\max_{M\times[0,t]}\big(\frac{|W|+|\nabla R|+|\nabla^2 R|}{R+b}+R+b\big).
\end{equation}
for all $(x,t)\in M\times[0,T)$.
\end{theorem}

At the end of this section, we finish the proof of Theorem \ref{thm1.4}.
\\\\\textbf{Proof of Theorem \ref{thm1.4}:} We prove this theorem by contradiction. Suppose that
\[\limsup_{t\rightarrow T}(\max_MR)<\infty\ \text{and}\ \limsup_{t\rightarrow T}\big(\max_M\frac{|W|+|\nabla R|+|\nabla^2 R|}{R}\big)<\infty.\]

Therefore, Theorem \ref{thm6.1} implies $|\overset\circ{Ric}|$ is uniformly bounded on $M\times[0,T)$. It is well known that the curvature on a $n$-dimensional $(n\geq3)$ Riemannian manifold can be decomposed into three orthogonal components the Weyl tensor part, the traceless Ricci part and the scalar curvature part (see e.g. \cite{Hsu19}). It follows that the Riemannian curvature is uniformly bounded on $M\times[0,T)$. By Theorem \ref{cor1.1}, this flow can be extended smoothly past time $T$, which contradicting the fact of $T$ is a maximal existence time.$\hfill\Box$
%\begin{remark}
%It is not hard to see from the calculation in Proposition \ref{prp2.1} and Lemma \ref{lem6.1} that if we assume $R_{min}(0)\geq0$ additionally in Theorem \ref{thm1.4},  the statement $\limsup_{t\rightarrow T}\big(\max_M\frac{|W|+|\nabla^2R|+|\nabla R|}{R}\big)=\infty$ can be
%\end{remark}
\section*{Acknowledgements}
The author would like to thank Prof. Yi Li and Prof. Xiaokui Yang for helpful suggestions.

\end{document}